\documentclass[12pt]{article}
\usepackage{amsmath}
\usepackage{epsf}
\usepackage{amssymb}
\usepackage{graphicx}
\usepackage{ifpdf}
\usepackage{caption}

\usepackage{array}
\usepackage{booktabs} 
\usepackage{multirow}

\newtheorem{theorem}{Theorem}[section]

\newtheorem{definition}{Definition}[section]
\newtheorem{lemma}{Lemma}[section]

\begin{document}
\pagestyle{empty}
\renewcommand{\thefootnote}{\fnsymbol{footnote}}

\begin{center}
{\bf \Large The generalized 4-connectivity of godan graphs}\footnote{This work was supported by National Science Foundation of China(No.12271157 \& No.12371346) and Natural Science Foundation of Hunan Province (No.2022JJ30028 \& No.2023JJ30072).}
\vskip 5mm

{{\bf Jing Wang$^1$, Yuanqiu Huang$^2$, Zhangdong Ouyang$^3$}\\[2mm]
$^1$ School of Mathematics, Changsha University, Changsha 410022, China\\
$^2$ School of Mathematics and Statistics, Hunan Normal University, Changsha 410081, China\\
$^3$ School of Mathematics and Statistics, Hunan First Normal University, Changsha 410205, China
}\\[6mm]
\end{center}
\date{}

\noindent{\bf Abstract}\; The generalized $k$-connectivity of a graph $G$, denoted by $\kappa_k(G)$, is the minimum number of internally edge disjoint $S$-trees for any $S\subseteq V(G)$ and $|S|=k$. The generalized $k$-connectivity is a natural extension of the classical connectivity and plays a key role in applications related to the modern interconnection networks. The godan graph $EA_n$ is a kind of Cayley graphs which posses many desirable properties. In this paper, we shall study the generalized 4-connectivity of $EA_n$ and show that $\kappa_4(EA_n)=n-1$ for $n\ge 3$.

{\bf MR(2000) Subject Classification} 05C40, 05C05

\section{Introduction}
\label{secintro}


A network is usually represented by a connected graph $G=(V(G),E(G))$, where $V(G)$ is the set of processors and $E(G)$ is the set of communication links between processors. Fault tolerance has become increasingly significant nowadays since multiprocessor systems failure is inevitable. It is important to consider how the network performs in the event of a certain number of nodes of failure and/or links failure in the network topology. The connectivity is an important parameter for measuring the fault tolerance of the network. A subset $S\subseteq V(G)$ of a connected graph $G$ is called a {\it vertex-cut} if $G\setminus S$ is disconnected or trivial. The {\it connectivity} $\kappa(G)$ of $G$ is defined as the minimum cardinality over all vertex-cuts of $G$. Note that the larger $\kappa(G)$ is, the more reliable the network is.

Let $P$ be a path in $G$ with $x$ and $y$ as its two terminal vertices, then $P$ is called an $(x,y)$-{\it path}. Two ($x,y$)-paths $P_1$ and $P_2$ are {\it internally disjoint} if they have no common internal vertices, that is, $V(P_1)\cap V(P_2)=\{x, y\}$. A well known theorem of Whitney \cite{Whitney1932} provides an equivalent definition of connectivity. For each 2-subset $S=\{x,y\}\subseteq V(G)$, let $\kappa(S)$ denote the maximum number of internally disjoint ($x,y$)-paths in $G$. Then
\begin{equation*}
\kappa(G)=\min\{\kappa(S) |S\subseteq V(G)\; {\rm and} \; |S|=2\}.
\end{equation*}

The generalized $k$-connectivity, which was introduced by Chartrand et al. \cite{Chartrand1984}, is a strengthening of connectivity and can be served as an essential parameter for measuring reliability and fault tolerance of the network. Let $G=(V(G),E(G))$ be a simple graph, $S$ be a subset of $V(G)$. A tree $T$ in $G$ is called an $S$-{\it tree}, if $S\subseteq V(T)$. The trees $T_1, T_2, \cdots, T_r$ are called {\it internally edge disjoint $S$-trees} if $V(T_i)\cap V(T_j)=S$ and $E(T_i)\cap E(T_j)=\emptyset$ for any integers $1\le i< j\le r$. Throughout this paper, we will use IDSTs to denote internally edge disjoint $S$-trees for short. Let $\kappa_G(S)$ denote the maximum number of IDSTs.

For an integer $k$ with $2\le k\le |V(G)|$, the {\it generalized $k$-connectivity} of $G$, denoted by $\kappa_k(G)$, is defined as
\begin{equation*}
\kappa_k(G)=\min\{\kappa_G(S) |S\subseteq V(G)\; {\rm and} \; |S|=k\}.
\end{equation*}

The generalized 2-connectivity is exactly the classical connectivity. Over the past few years, research on the generalized connectivity has received meaningful progress. Li et al. \cite{SLi2012n} derived that it is NP-complete for a general graph $G$ to decide whether there are $l$ internally edge disjoint trees connecting $S$, where $l$ is a fixed integer and $S\subseteq V(G)$. Authors in \cite{HZLi2014,SLi2010} investigated the upper and lower bounds of the generalized connectivity of a general graph $G$.

For $3\le k\le |V(G)|$, Li gave an upper bound of $\kappa_k(G)$ for a general graph $G$ in her Ph.D. thesis \cite{LiSSPhD2012}.

\begin{lemma}\label{lemupperKk}(\cite{LiSSPhD2012})
Let $G$ be a connected graph with minimum degree $\delta(G)$. If there are two adjacent vertices of degree $\delta(G)$, then $\kappa_k(G)\le \delta(G)-1$ for $3\le k\le |V(G)|$.
\end{lemma}

The following result is about the relationship between $\kappa_k(G)$ and $\kappa_{k-1}(G)$ of a regular graph $G$.

\begin{lemma}\label{lemKkk-1}(\cite{SLin2017})
Let $G$ be an $r$-regular graph. If $\kappa_k(G)=r-1$, then $\kappa_{k-1}(G)=r-1$, where $k\ge 4$.
\end{lemma}

Many authors tried to study exact values of the generalized connectivity of graphs. The generalized $k$-connectivity of the complete graphs, $\kappa_k(K_n)$, was determined in \cite{Chartrand2010} for every pair $k,n$ of integers with $2\le k\le n$. The generalized $k$-connectivity of the complete bipartite graphs $K_{a,b}$ was obtained in \cite{SLi2012b} for all $2\le k\le a+b$. Apart from these results, the generalized $k$-connectivity of hypercubes \cite{HZLi2012,SLin2017},  the generalized hypercubes \cite{Wang2021}, balanced hypercubes \cite{Wei2021}, exchanged hypercubes \cite{ZhaoHao20192}, dual cubes \cite{ZhaoHao20191}, star graphs and bubble-sort graphs \cite{SLi2016}, ($n,k$)-star graphs \cite{Snk2020,ANn12018},  bubble-sort star graphs \cite{Hao20191}, pancake graphs \cite{K4Pn}, Cayley graphs generated by trees and cycles \cite{SLi2017}, et al, have drawn many scholars' attention. So far, as we can see, the results on the generalized $k$-connectivity of networks were almost about $k\le 4$.

Ren and Wang \cite{godangraph2022} introduced the $n$-dimensional godan graph $EA_n$ as a topology structure of interconnection networks, they proved that $EA_n$ poses many desirable properties, such as regularity, high connectivity and vertex-transitivity. The authors in \cite{godangraph2022} also gave various connectivity, such as nature connectivity, 2-good-neighbour connectivity,  2-extra connectivity, of $EA_n$. The component connectivity of $EA_n$ was studied in \cite{GodanIn2022}. In this paper, we try to calculate the generalized 4-connectivity of $EA_n$ and obtain the following main result.

\begin{theorem}\label{thk4EAn}
Let $EA_n$ be the $n$-dimensional godan graph for $n\ge 3$. Then $\kappa_4(EA_n)=n-1$.
\end{theorem}

We end this section with some necessary preliminaries. We only consider a simple, connected graph $G=(V(G),E(G))$ with vertex set $V(G)$ and edge set $E(G)$. For $x\in V(G)$, the {\it degree} of $x$ in $G$, denoted by ${\rm deg}_G(x)$, is the number of edges of $G$ incident with $x$. Denote $\delta(G)$ the {\it minimum degree} of vertices of $G$. A graph is $d$-{\it regular} if ${\rm deg}_G(x)=d$ for every $x\in V(G)$. For $x\in V(G)$, we use $N_G(x)$ to denote the set of neighbours of $x$ in $G$. Let $V'\subseteq V(G)$, denote by $G\backslash V'$ the graph obtained from $G$ by deleting all the vertices in $V'$ together with their incident edges. Denote by $G[V']$ the subgraph of $G$ induced on $V'$.

The following Lemma \ref{lemxypath}, Lemma \ref{lemKfan} and Lemma \ref{lemXYpaths} are results on the connectivity of a graph that are well-known in the literature.

\begin{lemma}\label{lemxypath} (\cite{Bondy})
Let $G$ be a $k$-connected graph, and let $x$ and $y$ be a pair of distinct vertices of $G$. Then there exist $k$ internally disjoint ($x,y$)-paths in $G$.
\end{lemma}

\begin{lemma}\label{lemKfan} (\cite{Bondy})
Let $G$ be a $k$-connected graph, let $x$ be a vertex of $G$ and let $Y\subseteq V(G)\backslash\{x\}$ be a set of at least $k$ vertices of $G$. Then there exists a $k$-fan in $G$ from $x$ to $Y$, that is, there exists a family of $k$ internally disjoint ($x,Y$)-paths whose terminal vertices are distinct in $Y$.
\end{lemma}

\begin{lemma}\label{lemXYpaths} (\cite{Bondy})
Let $G$ be a $k$-connected graph, and let $X$ and $Y$ be subsets of $V(G)$ of cardinality at least $k$. Then there exists a family of $k$ pairwise disjoint ($X,Y$)-paths in $G$.
\end{lemma}

This paper is organized as follows. Section \ref{secdefEAn} presents the definition and properties of godan graphs. We note that the lower bound of $\kappa_4(EA_n)$ is critical to obtain the main result. Therefore, some lemmas which will be used to obtain the lower bound of $\kappa_4(EA_n)$ are given in Section \ref{seclemmas}. Then Theorem \ref{thk4EAn} is proved in Section \ref{secK4EAn}. Finally, the paper is concluded in Section \ref{seccon}.

\section{Definition and properties of godan graph}\label{secdefEAn}

Let $[n]=\{1,2,\cdots, n\}$. For convenience, we denote the permutation $\scriptsize{\begin{pmatrix} 1&2&\cdots & n\\ p_1&p_2&\cdots & p_n \end{pmatrix}}$ by $p_1p_2\cdots p_n$ and the permutation $\scriptsize{\begin{pmatrix} 1&\cdots& i&\cdots &j \cdots& n\\ 1&\cdots& j&\cdots &i \cdots& n\end{pmatrix}}$ by $(ij)$, while the latter is called a transposition. The composition $\sigma\circ\tau$ of two permutations $\sigma$ and $\tau$ is the function that maps any element $i$ to $\sigma(\tau(i))$.

Let $S_n$ be the symmetric group on $[n]$. The alternating group $A_n$ ($n\ge 3$) is the subgroup of $S_n$ containing all even permutations. Let $\Omega=\{(123),(132)\}\cup \{(12)(3i)| 4\le i\le n\}$. For $n\ge 3$, the $n$-dimensional alternating group network $AN_n$ (\cite{Ji1998}) is a graph with vertex set $V(AN_n)=A_n$ and edge set $E(AN_n)=\{uv| u=v\circ s, \;{\rm where}\; \{u,v\}\subseteq V(AN_n) \;{\rm and}\; s\in \Omega\}$.

\begin{figure}[htbp]
\begin{minipage}[t]{0.4\linewidth}
\centering
\resizebox{0.6\textwidth}{!} {\includegraphics{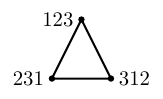}}
\caption{\small The alternating group network $AN_3$} \label{figAN3}
\end{minipage}
\begin{minipage}[t]{0.6\linewidth}
\centering
\resizebox{0.65\textwidth}{!} {\includegraphics{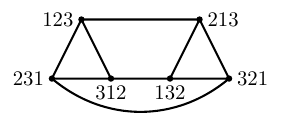}}
\caption{\small The godan graph $EA_3$} \label{figEA3}
\end{minipage}
\end{figure}

The alternating group network $AN_3$ is depicted in Figure \ref{figAN3}. It is worth mentioning that  $AN_n$ is in fact isomorphic to the $(n,n-2)$-star $S_{n,n-2}$ \cite{Cheng2012ANn}, whose generalized 4-connectivity is proved to be ($n-2$) in \cite{Snk2020}.

\begin{lemma}\label{lemANn}(\cite{Snk2020,ANn12019})
For the $n$-dimensional alternating group network $AN_n$ ($n\ge 3$), the following properties hold:\\
(1) $AN_n$ is ($n-1$)-regular with $\frac{n!}{2}$ vertices and $\frac{n!(n-1)}{4}$ edges;\\
(2) $\kappa(AN_n)=n-1$;\\
(3) $\kappa_4(AN_n)=\kappa_4(S_{n,n-2})=n-2$.
\end{lemma}

In \cite{godangraph2022}, Ren and Wang introduced the $n$-dimensional godan graph $EA_n$ as a topology structure of interconnection networks.

\begin{definition}\label{defEAn}(\cite{godangraph2022})
For $n\ge 3$, let $\Omega^*=\{(12),(123),(132)\}\cup \{(12)(3i)|4\le i\le n\}$. The $n$-dimensional godan graph $EA_n$ is the graph with vertex set $V(EA_n)=S_n$ in which two vertices $u$ and $v$ are adjacent in $EA_n$ if and only if $u=v\circ s^*$, where $s^*\in \Omega^*$.
\end{definition}

\begin{figure}[htbp]
\centering
\resizebox{0.7\textwidth}{!} {\includegraphics{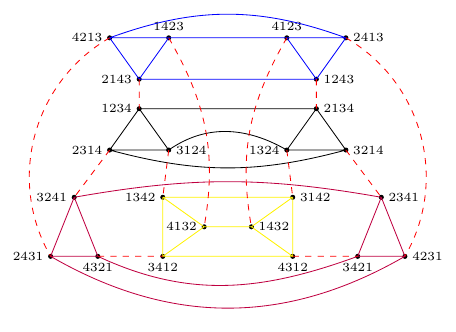}}
\caption{\small The godan graph $EA_4$} \label{figEA4}
\end{figure}

The godan graphs $EA_3$ and $EA_4$ are depicted in Figure \ref{figEA3} and Figure \ref{figEA4}, respectively. Some properties of $EA_n$ are listed in Lemma \ref{lemEAn1}.

\begin{lemma}\label{lemEAn1}(\cite{godangraph2022,GodanIn2022})
For $n\ge 3$, the $n$-dimensional godan graph $EA_n$ has the following properties:\\
(1) $EA_n$ is $n$-regular with $n!$ vertices and $\frac{n\times n!}{2}$ edges;\\
(2) $EA_n[A_n]\cong AN_n$, $EA_n[S_n\backslash A_n]\cong AN_n$;\\
(3) $E(A_n, S_n\backslash A_n)$ is a perfect matching of $EA_n$, where $E(A_n, S_n\backslash A_n)$ denotes all cross edges between $A_n$ and $S_n\backslash A_n$;\\
(4) $\kappa(EA_n)=n$.
\end{lemma}

It is seen that both $EA_n[A_n]$ and $EA_n[S_n\backslash A_n]$ are isomorphic to $AN_n$. For simplicity, we denote $EA_n[A_n]$ by $AN_n^1$ and denote $EA_n[S_n\backslash A_n]$ by $AN_n^2$, respectively. Then we can write the construction of $EA_n$ symbolically as
$$EA_n=AN_n^1\otimes AN_n^2.$$
For $1\le i\le 2$, let $x$ be a vertex in $AN_n^i$, it is clear that $x$ has a unique neighbor $x\circ (12)$ in $V(AN_n^{3-i})$. Denote $x\circ (12)$ by $\widetilde{x}$. We say that $\widetilde{x}$ is the {\it parity neighbour} of $x$. Due to Definition \ref{defEAn}, the following result is an easy observation.

\begin{lemma}\label{lemEAn2}
For $n\ge 3$, let $EA_n=AN_n^1\otimes AN_n^2$ and $x$ be a vertex of $AN_n^i$, $1\le i\le 2$. If $N_{AN_n^i}(x)=\{b_1, b_2, \cdots, b_{n-1}\}$. Then $N_{AN_n^{3-i}}(\widetilde{x})=\{\widetilde{b}_1, \widetilde{b}_2, \cdots, \widetilde{b}_{n-1}\}$.
\end{lemma}

\begin{lemma}\label{lemEAn31}
For $n\ge 4$, let $x$ be a vertex of $EA_n$ with $N_{EA_n}(x)=\{b_1, b_2, \cdots, b_{n-1}, \widetilde{x}\}$. Set $B=\{b_1, \cdots, b_{n-2}, \widetilde{x}, v\}$, where $v$ is a neighbour of $\widetilde{x}$ with $v\ne x$. Then $EA_n\backslash B$ is connected.
\end{lemma}

\noindent{\bf Proof}\; Remind that $EA_n=AN_n^1\otimes AN_n^2$. By symmetry, assume that $x\in V(AN_n^1)$. Let $B_1=\{b_1, \cdots, b_{n-2}\}$ and $B_2=\{\widetilde{x}, v\}$.

By Definition \ref{defEAn}, for $1\le i\le 2$, $B_i=B\cap V(AN_n^i)$ because $v\ne x$. 

By Lemma \ref{lemANn}(2), for $1\le i\le 2$, $AN_n^i\backslash B_i$ is connected since $|B_1|=n-2$ and $|B_2|=2$. Since $\frac{n!}{2}>n$ for $n\ge 4$ and $E(A_n, S_n\backslash A_n)$ is a perfect matching of $EA_n$, $EA_n\backslash B$ is still connected. \hfill$\Box$

\vskip 2mm
Now, we shall investigate the structure of $EA_n$ from another perspective. Let $EA_n^{m:i}$, called a {\it cluster}, be a subgraph of $EA_n$ induced on vertices with the $m$th symbol is $i$, where $4\le m\le n$ and $i\in[n]$. Note that each cluster $EA_n^{m:i}$ is isomorphic to $EA_{n-1}$, $4\le m\le n$ and $i\in[n]$. Thus, the construction of $EA_n$ may be described as
$$EA_n=EA_n^{m:1}\oplus EA_n^{m:2}\oplus \cdots \oplus EA_n^{m:n},\;\; 4\le m\le n.$$
For simplicity, we use $EA_n^{i}$ to denote $EA_n^{n:i}$ in the rest of this paper, $i\in[n]$. Thus, we may write the construction of $EA_n$ symbolically as
$$EA_n=EA_n^1\oplus EA_n^2\oplus \cdots \oplus EA_n^n.$$
For $i\in [n]$, let $x$ be a vertex in $EA_n^i$. Clearly, $x$ has ($n-1$) neighbours in $EA_n^i$ and exactly one neighbor $x'=x\circ(12)(3n)$ out of $EA_n^{i}$. We say that $x'$ is the {\it out neighbour} of $x$.

By Definition \ref{defEAn}, some properties of $EA_n$ can be derived directly, which are presented in Lemma \ref{lemEAn3}.

\begin{lemma}\label{lemEAn3}
For $n\ge 3$, let $EA_n=EA_n^1\oplus EA_n^2\oplus \cdots \oplus EA_n^n.$ The following properties hold:\\
(1) $|E(EA_n^i, EA_n^j)|=(n-2)!$,  where $i$ and $j$ are distinct integers in $[n]$.\\
(2) Let $x\in V(EA_n^i)$, $i\in [n]$. Then $|N_{EA_n^i}(x)|=n-1$. Furthermore, the vertices in $N_{EA_n^i}(x)$ can be ordered in an unique way $x_1, x_2, \cdots, x_{i-1}$, $x_{i+1}, \cdots, x_{n}$ such that $x_j'\in V(EA_n^j)$ for $j\in [n]\setminus \{i\}$. We use $P[x,x_j']$ to denote the path $\{xx_j, x_jx_j'\}$, where $j\in [n]\setminus \{i\}$.\\
(3) Let $x\in V(EA_n^i)$ and $x_j$ be a neighbour of $x$ in $EA_n^i$, $\{i,j\}\subseteq [n]$ and $j\ne i$. Then $x'$ and $x_j'$ belong to a same cluster of $EA_n$ if and only if $x_j=x\circ(12)$. Moreover, $x_j'=x'\circ(12)$.\\
(4) Let $x,y$ and $z$ be three distinct vertices of $EA_n$. The subgraph induced on $\{x,y,z\}$ is a 3-cycle if and only if $\{y,z\}=\{x\circ(123), x\circ(132)\}$.
\end{lemma}

We end this section with Lemma \ref{lemEAn4} which determines the lower bound of the connectivity of a subgraph of $EA_n$.

\begin{lemma}\label{lemEAn4}
For $n\ge 4$, let $EA_n=EA_n^1\oplus EA_n^2\oplus \cdots \oplus EA_n^n.$ Denote $H=\bigcup_{i=1}^{h}EA_n^{k_i}$ be a subgraph of $EA_n$, where $k_i\in [n]$ and $h\ge 2$. Then $\kappa(H)\ge n-2$.
\end{lemma}

\noindent{\bf Proof}\; Let $x$ and $y$ be any two distinct vertices in $H$. It suffices to prove that there are ($n-2$) internally disjoint ($x,y$)-paths in $H$.

{\bf Case 1}.  Both $x$ and $y$ belong to a same cluster, say $EA_n^{k_1}$.

By Lemma \ref{lemEAn1}(4), there are ($n-1$) internally disjoint ($x,y$)-paths in $EA_n^{k_1}$ since $\kappa(EA_n^{k_1})=\kappa(EA_{n-1})=n-1$.

{\bf Case 2}. $x$ and $y$ belong to different clusters of $EA_n$.

W.l.o.g., assume that $x\in V(EA_n^{k_1})$ and $y\in V(EA_n^{k_2})$. Based on Lemma \ref{lemEAn3}(1), $|E(EA_n^{k_1}, EA_n^{k_2})|=(n-2)!\ge n-2$ for $n\ge 4$. We may choose ($n-2$) different vertices $u_i\in V(EA_n^{k_1})$ such that $u_i'\in V(EA_n^{k_2})$, $i\in [n-2]$. Combined with Lemma \ref{lemKfan} and Lemma \ref{lemEAn1}(4), there is an ($n-2$)-fan $L_1, L_2, \cdots, L_{n-2}$ in $EA_n^{k_1}$ from $x$ to $\{u_1, u_2, \cdots, u_{n-2}\}$ such that $u_i\in V(L_i)$, where $i\in [n-2]$. Likewise, there is an ($n-2$)-fan $Q_1, Q_2, \cdots, Q_{n-2}$ in $EA_n^{k_2}$ from $y$ to $\{u_1', u_2', \cdots, u_{n-2}'\}$ such that $u_i'\in V(Q_i)$, $i\in [n-2]$.

For $i\in [n-2]$, let $P_i=L_i\cup Q_i\cup \{u_iu_i'\}$. Then $P_1, \cdots, P_{n-2}$ are ($n-2$) internally disjoint ($x,y$)-paths in $H$. \hfill$\Box$

\section{Some lemmas}\label{seclemmas}

The lower bound of $\kappa_4(EA_n)$ is the key to the proof of Theorem \ref{thk4EAn}. We shall make great efforts to provide several lemmas to obtain the lower bound of $\kappa_4(EA_n)$ for $n\ge 4$.

\begin{lemma}\label{lemS4}
For $n\ge 4$, let $EA_n=AN_n^1\otimes AN_n^2$ and $S=\{x,y,z,w\}$ be any 4-subset of $V(EA_n)$. If there exists an integer $i$ ($1\le i\le 2$) such that $S\subseteq V(AN_n^i)$. Then there are ($n-1$) IDSTs in $EA_n$.
\end{lemma}

\noindent{\bf Proof}\; W.l.o.g., assume that $S\subseteq V(AN_n^1)$. Then there are ($n-2$) IDSTs $T_1,\cdots, T_{n-2}$ in $AN_n^1$ by Lemma \ref{lemANn}(3) and Lemma \ref{lemEAn1}(2).

Note that $\{\widetilde{x}, \widetilde{y}, \widetilde{z}, \widetilde{w}\}\subseteq V(AN_n^2)$. There is a $\{\widetilde{x}, \widetilde{y}, \widetilde{z}, \widetilde{w}\}$-tree $T_{n-1}'$ in $AN_n^2$. Let $T_{n-1}=T_{n-1}'\cup \{x\widetilde{x}, y\widetilde{y}, z\widetilde{z}, w\widetilde{w}\}$. Clearly, $T_1,\cdots, T_{n-1}$ are ($n-1$) IDSTs in $EA_n$. \hfill$\Box$

\begin{lemma}\label{lemANS3}
For $n\ge 4$, let $EA_n=AN_n^1\otimes AN_n^2$ and $S=\{x,y,z,w\}$ be any 4-subset of $V(EA_n)$ satisfying the following three conditions:\\
(1) $\{x,y,z\}\subseteq V(AN_n^1)$ and $w\in V(AN_n^2)$;\\
(2) $\widetilde{w}\notin \{x,y,z\}$;\\
(3) $|E(\widetilde{w},\{x,y,z\})|\le 1$.\\
Then there are ($n-1$) IDSTs in $EA_n$.
\end{lemma}

\noindent{\bf Proof}\; Based on Lemma \ref{lemANn}(3), there are ($n-2$) internally edge disjoint $\{x,y,z,\widetilde{w}\}$-trees $T_{1}', \cdots, T_{n-2}'$ in $AN_n^1$. Furthermore, it has $1\le {\rm deg}_{T_i'}(\widetilde{w})\le 2$ for each $i\in [n-2]$ since ${\rm deg}_{AN_n^1}(\widetilde{w})=n-1$. For convenience, we may assume that ${\rm deg}_{T_i'}(\widetilde{w})=1$ and $\widetilde{u}_i$ be a neighbour of $\widetilde{w}$ in $T_i'$, where $i\in [n-3]$.

Let $\widetilde{U}=\{\widetilde{u}_1,\cdots, \widetilde{u}_{n-3}\}$ and $U=\{u_i|~ \widetilde{u}_i\in \widetilde{U}, ~i\in [n-3]\}.$

Note that $U\subseteq V(AN_n^2)$ and $|U|=n-3$.  Thus $AN_n^2\backslash U$ is still connected since $\kappa(AN_n^2)=n-1$. The following two cases are considered.

\vskip 2mm
{\bf Case 1}. $E(\widetilde{w}, \{x,y,z\})=\emptyset$.

Then $\{\widetilde{x},\widetilde{y},\widetilde{z}\}\cap \{w, u_1, \cdots, u_{n-3}\}=\emptyset$. Moreover, $wu_i\in E(AN_n^2)$ for $i\in [n-3]$ according to Lemma \ref{lemEAn2}.

There is a $\{\widetilde{x},\widetilde{y},\widetilde{z},w\}$-tree $T_{n-1}'$ in $AN_n^2\backslash U$. Let $T_{n-1}=T_{n-1}'\cup \{x\widetilde{x}, y\widetilde{y}, z\widetilde{z}\}$ and $T_{n-2}=T_{n-2}'\cup \{w\widetilde{w}\}$. For $i\in [n-3]$, let $T_{i}=(T_{i}'\backslash\{\widetilde{w}\})\cup \{wu_i, u_i\widetilde{u}_i\}$.

Then $T_{1}, \cdots, T_{n-1}$ are ($n-1$) IDSTs in $EA_n$.

\vskip 2mm
{\bf Case 2}. $|E(\widetilde{w}, \{x,y,z\})|=1$.

Assume that $x\widetilde{w}\in E(AN_n^1)$. If $x\notin \widetilde{U}$, then we can construct ($n-1$) IDSTs by similar methods in Case 1.

Hereafter, we consider the case that $x\in \widetilde{U}$. For simplicity, set $x=\widetilde{u}_1$. We only discuss the case that ${\rm deg}_{T_{n-2}'}(\widetilde{w})=2$, since ($n-1$) IDSTs can be obtained just by replacing $T_1$ and $T_{n-2}$ in Case 1 if ${\rm deg}_{T_{n-2}'}(\widetilde{w})=1$.

Let $\widetilde{u}_{n-2}$ and $\widetilde{u}_{n-1}$ be two neighbours of $\widetilde{w}$ in $T_{n-2}'$. Moreover, let $\widetilde{U}_1=\{\widetilde{u}_2,\cdots, \widetilde{u}_{n-1}\}$ and $U_1=\{u_i~|~ \widetilde{u}_i\in \widetilde{U}_1,  ~2\le i\le n-1\}.$ Then $U_1\subseteq V(AN_n^2)$ and $|U_1|=n-2$. Therefore, $AN_n^2\backslash U_1$ is connected and there is a $\{\widetilde{x},\widetilde{y},\widetilde{z},w\}$-tree $T_{n-1}'$ in $AN_n^2\backslash U_1$.

Let $$T_{n-1}=T_{n-1}'\cup \{x\widetilde{x}, y\widetilde{y}, z\widetilde{z}\},$$
$$T_{n-2}=(T_{n-2}'\backslash\{\widetilde{w}\})\cup \{wu_{n-2}, u_{n-2}\widetilde{u}_{n-2}, wu_{n-1}, u_{n-1}\widetilde{u}_{n-1}\}$$
and $T_{1}=T_{1}'\cup \{w\widetilde{w}\}$. For $2\le i\le n-3$, let $T_{i}$ be the same as in Case 1. Then $T_{1}, \cdots, T_{n-1}$ are ($n-1$) IDSTs in $EA_n$. \hfill$\Box$

\vskip 2mm

For $n\ge 4$ and $4\le m\le n$, recall that $EA_n=EA_n^{m:1}\oplus EA_n^{m:2}\oplus \cdots \oplus EA_n^{m:n}$. Let $S=\{x,y,z,w\}$ be any 4-subset of $V(EA_n)$. In order to obtain the lower bound of $\kappa_4(EA_n)$, the following subsections are distinguished based on the value of $\max\limits_{i\in [n]}\{|S\cap V(EA_n^{m:i})|\}$.

\subsection{$\max\limits_{i\in [n]}\{|S\cap V(EA_n^{m:i})|\}=3$, $4\le m\le n$.}

\begin{lemma}\label{lemS3}
For $n\ge 4$ and $4\le m\le n$, let $EA_n=EA_n^{m:1}\oplus EA_n^{m:2}\oplus \cdots \oplus EA_n^{m:n}$ and $S=\{x,y,z,w\}$ be any 4-subset of $V(EA_n)$. If there exists an integer $i\in[n]$ such that $|S\cap V(EA_n^{m:i})|=3$. Then there are ($n-1$) IDSTs in $EA_n$.
\end{lemma}

\noindent{\bf Proof}\; W.l.o.g., we assume that $\{x,y,z\}\subseteq V(EA_n^n)$ and $w\in V(EA_n^{n-1})$.

By Lemma \ref{lemEAn3}(2), for $i\in [n-2]$, there is a $\{x_i', y_i', z_i', w_i'\}$-tree $T_i'$ in $EA_n^i$ since $EA_n^i$ is connected, let
$$T_{i}=T_{i}'\cup P[x,x_{i}']\cup P[y,y_{i}']\cup P[z,z_{i}']\cup P[w,w_i'].$$
It is possible that $x_i=y_i$ for some $i\in [n-2]$. This possibility will not affect the result of this lemma.

Let $W=\{w_1, \cdots, w_{n-2}\}$. By Lemma \ref{lemEAn1}(4), $EA_n^{n-1}\backslash W$ is connected since $|W|=n-2$ and $\kappa(EA_n^{n-1})=n-1$. Note that there is a $\{x_{n-1}', y_{n-1}', z_{n-1}', w\}$-tree $T_{n-1}'$ in $EA_n^{n-1}\backslash W$. Let
$$T_{n-1}=T_{n-1}'\cup P[x,x_{n-1}']\cup P[y,y_{n-1}']\cup P[z,z_{n-1}'].$$

Denote $F$ the subgraph of $EA_n^n$ induced on $\{x,y,z\}$. Clearly, $|E(F)|\le 3$. The following cases need to be considered.

\vskip 2mm
{\bf Case 1}. $|E(F)|=0$.

Clearly, $T_1,\cdots, T_{n-1}$ are ($n-1$) desired IDSTs in $EA_n$.

\vskip 2mm
{\bf Case 2}. $|E(F)|=1$.

For simplicity, assume that $E(F)=\{xy\}$. Consider the out-neighbours of $x$ and $y$. The following discussions are based on whether $x'$ and $y'$ belong to a same cluster of $EA_n$ or not.

\vskip 2mm
{\bf Subcase 2.1}. There exists an integer $i\in[n-1]$ with $\{x',y'\}\subseteq V(EA_n^{i})$.

We discuss the case that $i=n-2$. For the other cases, ($n-1$) IDSTs can be obtained by similar arguments.

\begin{figure}[htbp]
\begin{minipage}[t]{0.5\linewidth}
\centering
\resizebox{0.95\textwidth}{!} {\includegraphics{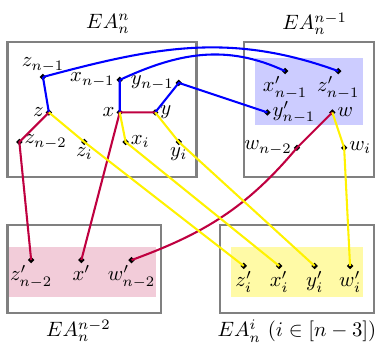}}
\caption{\small Illustration for Subcase 2.1} \label{figS3C21}
\end{minipage}
\begin{minipage}[t]{0.5\linewidth}
\centering
\resizebox{0.95\textwidth}{!} {\includegraphics{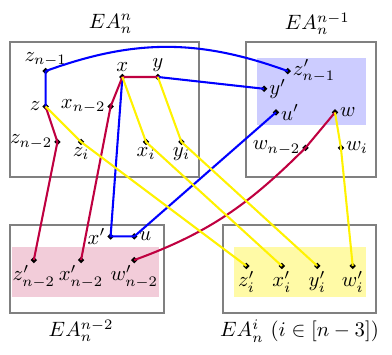}}
\caption{\small Illustration for Subcase 2.2} \label{figS3C22}
\end{minipage}
\end{figure}

It is inferred that $x=y_{n-2}$ and $y=x_{n-2}$. There is a $\{x', z_{n-2}', w_{n-2}'\}$-tree $T_{n-2}''$ in $EA_n^{n-2}$. Let
$$\widehat{T}_{n-2}=T_{n-2}''\cup P[z,z_{n-2}']\cup P[w,w_{n-2}']\cup\{xx', xy\}.$$
See Figure \ref{figS3C21}. Then $T_1,\cdots, T_{n-3}, \widehat{T}_{n-2}, T_{n-1}$ are ($n-1$) IDSTs in $EA_n$.

\vskip 2mm
{\bf Subcase 2.2}. $x'$ and $y'$ belong to different clusters of $EA_n$.

By symmetry, assume that $y'\in V(EA_n^{n-1})$ and $x'\in V(EA_n^{n-2})$. That means $x=y_{n-2}$ and $y=x_{n-1}$.

Recall that $x_{n-2}'\in V(EA_n^{n-2})$. Based on Lemma \ref{lemEAn3}(3) and (4), either $x'w_{n-2}'\notin E(EA_n^{n-2})$ or $x_{n-2}'w_{n-2}'\notin E(EA_n^{n-2})$. We only consider the case that $x'w_{n-2}'\notin E(EA_n^{n-2})$. For the other case, the proof is similar. By Lemma \ref{lemEAn3}(2), there is a neighbour of $x'$ in $EA_n^{n-2}$, say $u$, such that $u'\in V(EA_n^{n-1})$. Clearly, $u\notin \{x_{n-2}',z_{n-2}',w_{n-2}'\}$.

There exists a $\{u',y',z_{n-1}',w\}$-tree $T_{n-1}''$ in $EA_n^{n-1}\backslash W$. Meanwhile, there is a $\{x_{n-2}',z_{n-2}',w_{n-2}'\}$-tree $T_{n-2}''$ in $EA_n^{n-2}\backslash \{x',u\}$ since $EA_n^{n-2}\backslash \{x',u\}$ is connected.

Let
$$\widehat{T}_{n-2}=T_{n-2}''\cup P[x,x_{n-2}']\cup P[z,z_{n-2}']\cup P[w,w_{n-2}']\cup \{xy\}$$
and
$$\widehat{T}_{n-1}=T_{n-1}''\cup P[z,z_{n-1}']\cup \{xx', x'u, uu', yy'\}.$$
See Figure \ref{figS3C22}. Then $T_1,\cdots, T_{n-3}, \widehat{T}_{n-2}, \widehat{T}_{n-1}$ are ($n-1$) IDSTs in $EA_n$.

\vskip 2mm
{\bf Case 3}. $|E(F)|=2$.

For simplicity, assume that $E(F)=\{xy, xz\}$. By Lemma \ref{lemEAn3}(2), $y'$ and $z'$ belong to different clusters of $EA_n$ since both $y$ and $z$ are neighbours of $x$ in $EA_n^n$.

\vskip 2mm
{\bf Subcase 3.1}. $\{y', z'\}\cap V(EA_n^{n-1})\ne \emptyset$.

Assume that $z'\in V(EA_n^{n-1})$ and $y'\in V(EA_n^{n-2})$. That means $z=x_{n-1}$ and $y=x_{n-2}$ . Consider now the location of $x'$.

\begin{figure}[htbp]
\begin{minipage}[t]{0.5\linewidth}
\centering
\resizebox{0.95\textwidth}{!} {\includegraphics{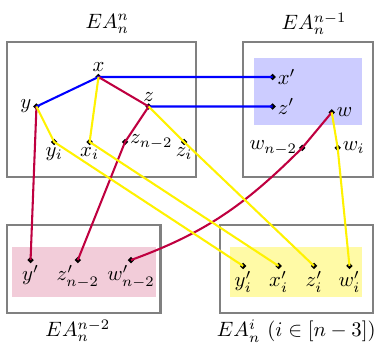}}
\caption{\small Illustration for Subcase 3.1.1} \label{figS3C311}
\end{minipage}
\begin{minipage}[t]{0.5\linewidth}
\centering
\resizebox{0.95\textwidth}{!} {\includegraphics{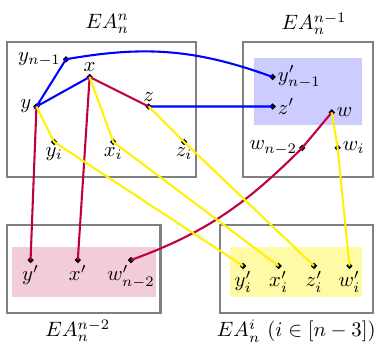}}
\caption{\small Illustration for Subcase 3.1.2} \label{figS3C312}
\end{minipage}
\end{figure}

{\bf Subcase 3.1.1}. $x'\in V(EA_n^{n-1})$.

Then $x=y_{n-1}=z_{n-1}$. There is a $\{x',z',w\}$-tree $T_{n-1}''$ in $EA_n^{n-1}\backslash W$ and a $\{y', z_{n-2}', w_{n-2}'\}$-tree $T_{n-2}''$ in $EA_n^{n-2}$, respectively. Let
$$\widehat{T}_{n-2}=T_{n-2}''\cup P[z,z_{n-2}']\cup P[w,w_{n-2}']\cup \{xz, yy'\}$$
and
$$\widehat{T}_{n-1}=T_{n-1}''\cup \{xy, xx', zz'\}.$$
See Figure \ref{figS3C311}. It is seen that $T_1,\cdots, T_{n-3}, \widehat{T}_{n-2}, \widehat{T}_{n-1}$ are ($n-1$) IDSTs in $EA_n$.

\vskip 2mm
{\bf Subcase 3.1.2}. $x'\in V(EA_n^{n-2})$.

That means $x=y_{n-2}=z_{n-2}$. There is a $\{x',y',w_{n-2}'\}$-tree $T_{n-2}''$ in $EA_{n}^{n-2}$ and a $\{y_{n-1}', z', w\}$-tree $T_{n-1}''$ in $EA_{n}^{n-1}\backslash W$, respectively. Let
$$\widehat{T}_{n-2}=T_{n-2}''\cup P[w,w_{n-2}']\cup \{xx', xz, yy'\}$$
and
$$\widehat{T}_{n-1}=T_{n-1}''\cup P[y,y_{n-1}']\cup \{xy, zz'\}.$$
See Figure \ref{figS3C312}. Clearly, $T_1,\cdots, T_{n-3}, \widehat{T}_{n-2}, \widehat{T}_{n-1}$ are ($n-1$) IDSTs in $EA_n$.

\vskip 2mm
{\bf Subcase 3.1.3}. $x'\in \bigcup_{i=1}^{n-3}V(EA_n^{i})$.

W.l.o.g., assume that $x'\in V(EA_n^{n-3})$, which means that $x=y_{n-3}=z_{n-3}$.

Recall that $x_{n-3}'\in V(EA_n^{n-3})$. According to Lemma \ref{lemEAn3}(3) and (4), either $x'w_{n-3}'\notin E(EA_n^{n-3})$ or $x_{n-3}'w_{n-3}'\notin E(EA_n^{n-3})$. For convenience, we only consider the case that $x'w_{n-3}'\notin E(EA_n^{n-3})$. Let $u$ be a neighbour of $x'$ in $EA_n^{n-3}$ such that $u'\in V(EA_n^{n-1})$ and $v$ be a neighbour of $y'$ in $EA_n^{n-2}$ such that $v'\in V(EA_n^{n-3})$. It can be inferred that $u\ne w_{n-3}'$.

Let $F_1=EA_{n}^{n-2}\backslash \{y',v\}$ and $F_2=EA_{n}^{n-3}\backslash \{x',u\}.$

\begin{figure}[htbp]
\begin{minipage}[t]{0.5\linewidth}
\centering
\resizebox{0.95\textwidth}{!} {\includegraphics{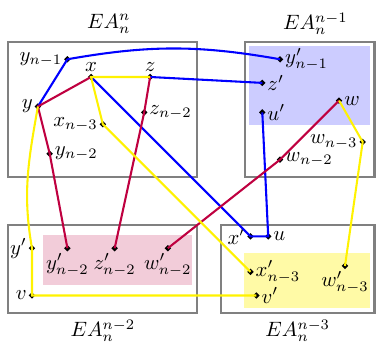}}
\caption{\small Illustration for Subcase 3.1.3} \label{figS3C313}
\end{minipage}
\begin{minipage}[t]{0.5\linewidth}
\centering
\resizebox{0.95\textwidth}{!} {\includegraphics{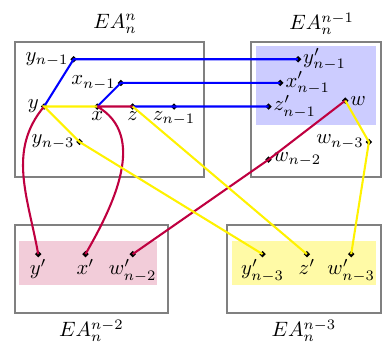}}
\caption{\small Illustration for Subcase 3.2.1} \label{figS3C321}
\end{minipage}
\end{figure}

Clearly, both $F_1$ and $F_2$ are connected. Thus, there is a $\{y_{n-2}',z_{n-2}',w_{n-2}'\}$-tree $T_{n-2}''$ in $F_1$ and a $\{x_{n-3}',v',w_{n-3}'\}$-tree $T_{n-3}''$ in $F_2$, respectively. Moreover, there is a $\{y_{n-1}',z',w,u'\}$-tree $T_{n-1}''$ in $EA_{n}^{n-1}\backslash W$. Let
$$\widehat{T}_{n-3}=T_{n-3}''\cup P[x,x_{n-3}']\cup P[w,w_{n-3}']\cup \{xz, yy', y'v, vv'\},$$
$$\widehat{T}_{n-2}=T_{n-2}''\cup P[y,y_{n-2}']\cup P[z,z_{n-2}']\cup P[w,w_{n-2}']\cup \{xy\}$$
and
$$\widehat{T}_{n-1}=T_{n-1}''\cup P[y,y_{n-1}']\cup\{xx',x'u,uu',zz'\}.$$
The trees $\widehat{T}_{n-3}, \widehat{T}_{n-2}$ and $\widehat{T}_{n-1}$ are depicted in Figure \ref{figS3C313}. Thus, $T_1,\cdots, T_{n-4}, \widehat{T}_{n-3}, \widehat{T}_{n-2},$ $\widehat{T}_{n-1}$ are ($n-1$) IDSTs in $EA_n$.

{\bf Subcase 3.2}. $\{y', z'\}\cap V(EA_n^{n-1})=\emptyset$.

For simplicity, let $y'\in V(EA_n^{n-2})$ and $z'\in V(EA_n^{n-3})$. The next discussions are based on the location of $x'$.

\vskip 2mm
{\bf Subcase 3.2.1}. Either $x'\in V(EA_{n}^{n-2})$ or $x'\in V(EA_{n}^{n-3})$.

By symmetry, assume that $x'\in V(EA_{n}^{n-2})$. There is a $\{x',y',w_{n-2}'\}$-tree $T_{n-2}''$ in $EA_{n}^{n-2}$ and a $\{y_{n-3}',z',w_{n-3}'\}$-tree $T_{n-3}''$ in $EA_{n}^{n-3}$, respectively. Let
$$\widehat{T}_{n-3}=T_{n-3}''\cup P[y,y_{n-3}']\cup P[w,w_{n-3}']\cup \{xy,zz'\}$$
and
$$\widehat{T}_{n-2}=T_{n-2}''\cup P[w,w_{n-2}']\cup \{xx', xz, yy'\}.$$
See Figure \ref{figS3C321}.  Thus, $T_1,\cdots, T_{n-4}, \widehat{T}_{n-3}, \widehat{T}_{n-2}, T_{n-1}$ are ($n-1$) IDSTs in $EA_n$.

\vskip 2mm
{\bf Subcase 3.2.2}. $x'\in V(EA_{n}^{i})$, where $i\in [n-4]$ or $i=n-1$.

We consider the case that $i=n-1$. For the other cases, the proofs are similar.

By similar arguments in Subcase 2.2 and Lemma \ref{lemEAn3}, either $y'w_{n-2}'\notin E(EA_{n}^{n-2})$ or $y_{n-2}'w_{n-2}'\notin E(EA_{n}^{n-2})$. For simplicity, assume that $y'w_{n-2}'\notin E(EA_{n}^{n-2})$. Likewise, we may assume that $z'w_{n-3}'\notin E(EA_{n}^{n-3})$.

Let $u$ be a neighbour of $y'$ in $EA_{n}^{n-2}$ such that $u'\in V(EA_{n}^{n-1})$ and $v$ be a neighbour of $z'$ in $EA_{n}^{n-3}$ such that $v'\in V(EA_{n}^{n-1})$, respectively. It has $u\ne w_{n-2}'$ and $v\ne w_{n-3}'$ by former analysis.

Let $F_1=EA_{n}^{n-2}\backslash \{y',u\}$ and $F_2=EA_{n}^{n-3}\backslash \{z',v\}.$

\begin{figure}[htbp]
\begin{minipage}[t]{0.5\linewidth}
\centering
\resizebox{0.95\textwidth}{!} {\includegraphics{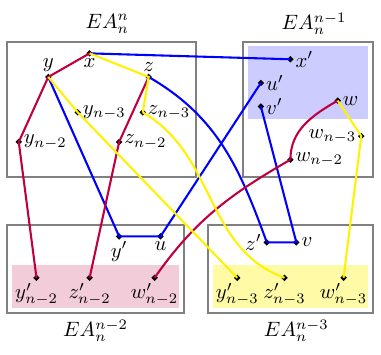}}
\caption{\small Illustration for Subcase 3.2.2} \label{figS3C322}
\end{minipage}
\begin{minipage}[t]{0.5\linewidth}
\centering
\resizebox{0.95\textwidth}{!} {\includegraphics{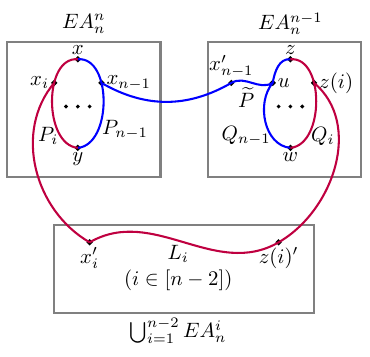}}
\caption{\small Illustration for Case 1} \label{figS22C1}
\end{minipage}
\end{figure}

We may obtain a $\{y_{n-2}',z_{n-2}',w_{n-2}'\}$-tree $T_{n-2}''$ in $F_1$ and a $\{y_{n-3}',z_{n-3}',w_{n-3}'\}$-tree $T_{n-3}''$ in $F_2$ since both $F_1$ and $F_2$ are connected. Moreover, there is a $\{x',u',v',w\}$-tree $T_{n-1}''$ in $EA_{n}^{n-1}\backslash W$. Let
$$\widehat{T}_{n-3}=T_{n-3}''\cup P[y,y_{n-3}']\cup P[z,z_{n-3}']\cup P[w,w_{n-3}']\cup \{xz\},$$
$$\widehat{T}_{n-2}=T_{n-2}''\cup P[y,y_{n-2}']\cup P[z,z_{n-2}']\cup P[w,w_{n-2}']\cup \{xy\}$$
and
$$\widehat{T}_{n-1}=T_{n-1}''\cup \{xx', yy', y'u, uu', zz', z'v, vv'\}.$$
See Figure \ref{figS3C322}. Therefore, $T_1,\cdots, T_{n-4}, \widehat{T}_{n-3}, \widehat{T}_{n-2}, \widehat{T}_{n-1}$ are ($n-1$) IDSTs in $EA_n$.

{\bf Case 4}. $|E(F)|=3$.

That is to say, $E(F)=\{xy, xz, yz\}$. Remind that $EA_n=AN_n^1\otimes AN_n^2$. Either $\{x,y,z\}\subseteq V(AN_n^1)$ or $\{x,y,z\}\subseteq V(AN_n^2)$ by Lemma \ref{lemEAn3}(4). By symmetry, assume that $\{x,y,z\}\subseteq V(AN_n^1)$. With the help of Lemma \ref{lemS4}, we only need to consider the possibility that $w\in V(AN_n^2)$.

Consider now $\widetilde{w}$, the parity neighbour of $w$. Claim 1 can be obtained easily.

\vskip 2mm
\noindent{\bf Claim 1}. $\widetilde{w}\notin \{x,y,z\}$.

\noindent{\bf Proof of Claim 1}. Suppose to contrary that $\widetilde{w}=x$. Then $w=\widetilde{x}\in V(EA_n^n)$. A contradiction to the assumption of this lemma.  \hfill$\Box$

\vskip 2mm
Based on Lemma \ref{lemEAn3}(4), Claim 2 may be followed immediately.

\vskip 2mm
\noindent{\bf Claim 2}. $|E(\widetilde{w}, \{x,y,z\})|\le 1$.

Therefore, Lemma \ref{lemANS3} implies that there are ($n-1$) IDSTs in $EA_n$.

In all, the proof is completed. \hfill$\Box$

\subsection{$\max\limits_{i\in [n]}\{|S\cap V(EA_n^{m:i})|\}=2$, $4\le m\le n$.}

\begin{lemma}\label{lemS22}
For $n\ge 4$ and $4\le m\le n$, let $EA_n=EA_n^{m:1}\oplus EA_n^{m:2}\oplus \cdots \oplus EA_n^{m:n}$ and $S=\{x,y,z,w\}$ be any 4-subset of $V(EA_n)$. If there exist different integers $i$ and $j$ in $[n]$ such that $|S\cap V(EA_n^{m:i})|=|S\cap V(EA_n^{m:j})|=2$. Then there are ($n-1$) IDSTs in $EA_n$.
\end{lemma}

\noindent{\bf Proof}\; W.l.o.g., assume that $\{x,y\}\subseteq V(EA_n^n)$ and $\{z,w\}\subseteq V(EA_n^{n-1})$.

By Lemma \ref{lemEAn1}(4), $\kappa(EA_n^i)=\kappa(EA_{n-1})=n-1$ for each $i\in [n]$. Combined with Lemma \ref{lemxypath}, there are ($n-1$) internally disjoint ($x,y$)-paths $P_1, \cdots, P_{n-1}$ in $EA_n^n$ and ($n-1$) internally disjoint ($z,w$)-paths $Q_1, \cdots, Q_{n-1}$ in $EA_n^{n-1}$, respectively. For convenience, let $x_i\in V(P_i)$ for $i\in [n-1]$.

Since $EA_n^{n-1}$ is connected, there is a ($x_{n-1}',z$)-path $\widetilde{P}$ in $EA_n^{n-1}$. Let $u$ be the first vertex in $V(\widetilde{P})\cap (\bigcup_{i=1}^{n-1}V(Q_i))$ when $\widetilde{P}$ starts from $x_{n-1}'$. For convenience, assume that $u\in V(Q_{n-1})$.

Let $z(i)$ be the neighbour of $z$ such that $z(i)\in V(Q_i)$, $i\in [n-1]$. Set
\begin{equation*}
X=\{x_1',\cdots, x_{n-2}'\} \;\;{\rm and}\;\; Z=\{z(1)',\cdots, z(n-2)'\}.
\end{equation*}
According to Lemma \ref{lemEAn3}(2), $X\cap V(EA_n^{n-1})=\emptyset$ and $|Z\cap V(EA_n^{n})|\le 1$.

\vskip 2mm
{\bf Case 1}. $Z\cap V(EA_n^{n})=\emptyset$.

\vskip 2mm
By Lemma \ref{lemEAn4}, $\kappa\big(\bigcup_{i=1}^{n-2}EA_n^i\big)\ge n-2$. With the help of Lemma \ref{lemXYpaths}, there are ($n-2$) internally disjoint ($X,Z$)-paths $L_1,\cdots, L_{n-2}$ in $\bigcup_{i=1}^{n-2}EA_n^i$ such that $\{x_i',z(i)'\}\subseteq V(L_i)$ for $i\in[n-2]$. For $i\in [n-2]$, let
$$T_i=P_i\cup Q_i\cup L_i\cup\{x_ix_i',z(i)z(i)'\}.$$
Let $T_{n-1}=P_{n-1}\cup Q_{n-1}\cup \widetilde{P}(x_{n-1}',u)\cup \{x_{n-1}x_{n-1}'\}$. See Figure \ref{figS22C1}. Then $T_1,\cdots, T_{n-1}$ are ($n-1$) IDSTs in $EA_n$.

\vskip 2mm
{\bf Case 2}. $|Z\cap V(EA_n^{n})|=1$.

For simplicity, we may assume that $z(n-2)'\in V(EA_n^{n})$.

\vskip 2mm
{\bf Subcase 2.1}. $z(n-2)'\in \bigcup_{i=1}^{n-1}V(P_i)$.

\vskip 2mm
{\bf Subcase 2.1.1}. $z(n-2)'\in V(P_{n-1})$.

Let $Z_1=\{z(1)',\cdots, z(n-3)',z(n-1)'\}$. Note that $Z_1\cap V(EA_n^{n})=\emptyset$. By similar arguments in Case 1, there are ($n-2$) internally disjoint ($X,Z_1$)-paths $L_1,\cdots, L_{n-2}$ in $\bigcup_{i=1}^{n-2}EA_n^i$ such that $\{x_i',z(i)'\}\subseteq V(L_i)$ for $i\in[n-3]$ and $\{x_{n-2}',z(n-1)'\}\subseteq V(L_{n-2})$.

Let $T_{n-2}=P_{n-2}\cup Q_{n-1}\cup L_{n-2}\cup \{x_{n-2}x_{n-2}', z(n-1)z(n-1)'\}$ and $T_{n-1}=P_{n-1}\cup Q_{n-2}\cup \{z(n-2)z(n-2)'\}$. For $i\in [n-3]$, we may construct $T_i$ be the same as in Case 1. Then $T_1,\cdots, T_{n-1}$ are ($n-1$) IDSTs in $EA_n$.

\vskip 2mm
{\bf Subcase 2.1.2}. $z(n-2)'\in \bigcup_{i=1}^{n-2}V(P_i)$.

For convenience, assume that $z(n-2)'\in V(P_{n-2})$. Let $T_{n-2}=P_{n-2}\cup Q_{n-2}\cup \{z(n-2)z(n-2)'\}$ and $T_i$ be the same as in Case 1 for $i\in [n-3]$ or $i=n-1$. Then $T_1,\cdots, T_{n-1}$ are ($n-1$) IDSTs in $EA_n$.

\vskip 2mm
{\bf Subcase 2.2}. $z(n-2)'\notin \bigcup_{i=1}^{n-1}V(P_i)$.

Let $u_1,\cdots, u_{n-1}$ be ($n-1$) neighbours of $z(n-2)'$ in $EA_n^n$. Firstly, we consider the case that $u_1\in \bigcup_{i=1}^{n-1}V(P_i)$. By similar analysis in Subcase 2.1, we can obtain ($n-1$) IDSTs in $EA_n$ by replacing the edge $z(n-2)z(n-2)'$ by the path $\{u_1z(n-2)', z(n-2)z(n-2)'\}$ in $T_{n-1}$ when $u_1\in V(P_{n-1})$ or in $T_{n-2}$ when $u_1\in V(P_{n-2})$.

Then we consider the case that $\{u_1,\cdots, u_{n-1}\}\cap \big(\bigcup_{i=1}^{n-1}V(P_i)\big)=\emptyset$. By Lemma \ref{lemEAn3}(2), there is a vertex in $\{u_1,\cdots, u_{n-1}\}$, say $u_1$, such that $u_1'\in \bigcup_{i=1}^{n-2}V(EA_n^i)$. Let
$$Z_2=\{z(1)',\cdots, z(n-3)',u_1'\}.$$
Clearly, $|Z_2|=n-2$ and $Z_2\cap V(EA_n^n)=\emptyset$.

Combined with Lemma \ref{lemXYpaths} and Lemma \ref{lemEAn4}, there are ($n-2$) internally disjoint ($X,Z_2$)-paths $L_1, \cdots, L_{n-2}$ in $\bigcup_{i=1}^{n-2}EA_n^i$ such that $\{x_i',z(i)'\}\subseteq V(L_i)$ for $i\in[n-3]$ and $\{x_{n-2}',u_1'\}\subseteq V(L_{n-2})$. Let
$$T_{n-2}=P_{n-2}\cup Q_{n-2}\cup L_{n-2}\cup \{x_{n-2}x_{n-2}', z(n-2)z(n-2)', z(n-2)'u_1, u_1u_1'\}$$
and $T_i$ be the same as in Case 1 for $i\in [n-3]$ or $i=n-1$. Then $T_1,\cdots, T_{n-1}$ are ($n-1$) IDSTs in $EA_n$. The proof is done. \hfill$\Box$

\begin{lemma}\label{lemS211}
For $n\ge 4$ and $4\le m\le n$, let $EA_n=EA_n^{m:1}\oplus EA_n^{m:2}\oplus \cdots \oplus EA_n^{m:n}$ and $S=\{x,y,z,w\}$ be any 4-subset of $V(EA_n)$. If there exist different integers $i$, $j$ and $k$ in $[n]$ such that $|S\cap V(EA_n^{m:i})|=2$ and $|S\cap V(EA_n^{m:j})|=|S\cap V(EA_n^{m:k})|=1$. Then there are ($n-1$) IDSTs in $EA_n$.
\end{lemma}

\noindent{\bf Proof}\; W.l.o.g., assume that $\{x,y\}\subseteq V(EA_n^n)$, $z\in V(EA_n^{n-1})$ and $w\in V(EA_n^{n-2})$.

There are ($n-1$) internally disjoint ($x,y$)-paths $P_1, \cdots, P_{n-1}$ in $EA_n^n$. For simplicity, let $x_i\in V(P_i)$, $i\in[n-1]$.

For $i\in[n-3]$, there is a $\{x_i',z_i',w_i'\}$-tree $T_{i}'$ in $EA_n^i$ since $\{x_i',z_i',w_i'\}\subseteq V(EA_n^i)$. Therefore, let $T_i=T_{i}'\cup P_i\cup P[z,z_i']\cup P[w,w_i']\cup \{x_ix_i'\}$ for $i\in[n-3]$.

Let $Z=\{z_1, \cdots, z_{n-2}\}$ and $W=\{w_1, \cdots, w_{n-3}, w_{n-1}\}$.

Clearly, $|Z|=|W|=n-2$. Both $EA_n^{n-1}\backslash Z$ and $EA_n^{n-2}\backslash W$ are connected since $\kappa(EA_n^{i})=n-1$ for $i\in [n]$.

The following discussions are based on whether $w_{n-1}=z_{n-2}'$ or not.

\vskip 2mm
{\bf Case 1}. $w_{n-1}\ne z_{n-2}'$.

\begin{figure}[htbp]
\begin{minipage}[t]{0.5\linewidth}
\centering
\resizebox{0.95\textwidth}{!} {\includegraphics{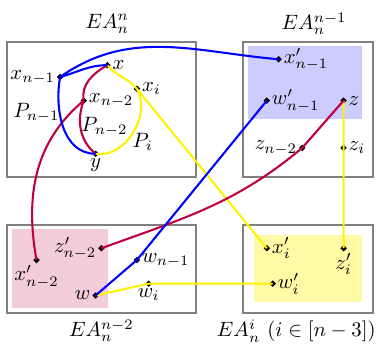}}
\caption{\small Illustration for Case 1} \label{figS211C1}
\end{minipage}
\begin{minipage}[t]{0.5\linewidth}
\centering
\resizebox{0.93\textwidth}{!} {\includegraphics{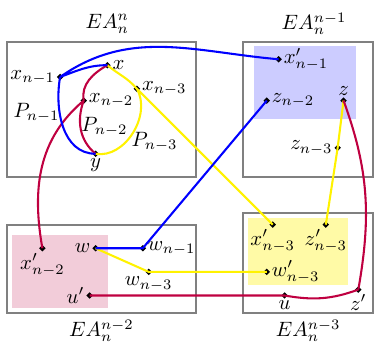}}
\caption{\small Illustration for Subcase 2.1} \label{figS211C21}
\end{minipage}
\end{figure}

There is a $\{x_{n-1}',z,w_{n-1}'\}$-tree $T_{n-1}'$ in $EA_n^{n-1}\backslash Z$ and a $\{x_{n-2}',z_{n-2}',w\}$-tree $T_{n-2}'$ in $EA_n^{n-2}\backslash W$, respectively. Let
$$T_{n-1}=T_{n-1}'\cup P_{n-1}\cup P[w,w_{n-1}']\cup \{x_{n-1}x_{n-1}'\}$$
and
$$T_{n-2}=T_{n-2}'\cup P_{n-2}\cup P[z,z_{n-2}']\cup \{x_{n-2}x_{n-2}'\}.$$
See Figure \ref{figS211C1}. Then $T_1,\cdots, T_{n-1}$ are ($n-1$) IDSTs in $EA_n$.

\vskip 2mm
{\bf Case 2}. $w_{n-1}=z_{n-2}'$.

Now consider the location of $z'$.

\vskip 2mm
{\bf Subcase 2.1}. $z'\in \bigcup_{i=1}^{n-3}V(EA_n^i)$.

W.l.o.g., assume that $z'\in V(EA_n^{n-3})$. By similar analysis in Subcase 2.2 of Lemma \ref{lemS3}, either $z'w_{n-3}'\notin E(EA_n^{n-3})$ or $z_{n-3}'w_{n-3}'\notin E(EA_n^{n-3})$. We only consider the possibility that $z'w_{n-3}'\notin E(EA_n^{n-3})$ for convenience. Thus, there is a neighbour of $z'$ in $EA_n^{n-3}$, say $u$, such that $u'\in V(EA_n^{n-2})$.

There is a $\{x_{n-1}', z, w_{n-1}'\}$-tree $T_{n-1}'$ in $EA_n^{n-1}\backslash \{z_1, \cdots, z_{n-3}\}$, a $\{x_{n-2}', w, u'\}$-tree $T_{n-2}'$ in $EA_n^{n-2}\backslash W$ and a $\{x_{n-3}', z_{n-3}', w_{n-3}'\}$-tree $T_{n-3}''$ in $EA_n^{n-3}\backslash \{z', u\}$, respectively.

Let
$$T_{n-1}=T_{n-1}'\cup P_{n-1}\cup P[w,w_{n-1}']\cup \{x_{n-1}x_{n-1}'\},$$
$$T_{n-2}=T_{n-2}'\cup P_{n-2}\cup \{x_{n-2}x_{n-2}', zz', z'u, uu'\}$$
and
$$\widehat{T}_{n-3}=T_{n-3}''\cup P_{n-3}\cup P[z,z_{n-3}']\cup P[w,w_{n-3}']\cup \{x_{n-3}x_{n-3}'\}.$$
See Figure \ref{figS211C21}. Then $T_1, \cdots, T_{n-4},$ $\widehat{T}_{n-3}, T_{n-2}, T_{n-1}$ are ($n-1$) IDSTs in $EA_n$.

\vskip 2mm
{\bf Subcase 2.2}. $z'\in V(EA_n^{n-2})$.

There is a $\{x_{n-2}', z', w\}$-tree $T_{n-2}'$ in $EA_n^{n-2}\backslash W$. Let $T_{n-2}=T_{n-2}'\cup P_{n-2}\cup \{x_{n-2}x_{n-2}', zz'\}$. For $i\in [n-3]$ or $i=n-1$, let $T_i'$ and $T_i$ be the same as in Case 1. Then $T_1, \cdots,  T_{n-1}$ are ($n-1$) IDSTs in $EA_n$.

\vskip 2mm
{\bf Subcase 2.3}. $z'\in V(EA_n^{n})$.

\vskip 2mm
{\bf Subcase 2.3.1}. $z'\in \bigcup_{i=1}^{n-1}V(P_i)$.

By Lemma \ref{lemEAn3}(3), it has $z_n=z\circ(12)$. There is a neighbour of $z_n$ in $EA_n^{n-1}$, say $v$, such that $v'\in V(EA_n^{n-2})$. According to Lemma \ref{lemEAn3}(4), $v\ne z_{n-2}$.

\vskip 2mm
{\bf Subcase 2.3.1.1}. $z'\in V(P_{n-1})$.

There is a $\{x_{n-2}', v', w\}$-tree $T_{n-2}'$ in $EA_n^{n-2}\backslash W$. Let $T_{n-2}=T_{n-2}'\cup P_{n-2}\cup \{x_{n-2}x_{n-2}', zz_n, z_nv, vv'\}$. Let $T_{n-1}=P_{n-1}\cup P[w,w_{n-1}']\cup \{zz', zw_{n-1}'\}$ since $w_{n-1}'=z_{n-2}$. Then $T_1, \cdots,  T_{n-1}$ are ($n-1$) IDSTs in $EA_n$.

\vskip 2mm
{\bf Subcase 2.3.1.2}. $z'\in V(P_{n-2})$.

There is a $\{x_{n-1}', z, z_n, v\}$-tree $T_{n-1}'$ in $EA_{n}^{n-1}\backslash Z$ and a ($w, v'$)-path $\widetilde{P}$ in $EA_{n}^{n-2}\backslash W$. Let $T_{n-2}=P_{n-2}\cup P[w,w_{n-1}']\cup \{zz', zw_{n-1}'\}$ and $T_{n-1}=T_{n-1}'\cup P_{n-1}\cup \widetilde{P}\cup \{x_{n-1}x_{n-1}', vv'\}$. Then $T_1, \cdots,  T_{n-1}$ are ($n-1$) IDSTs in $EA_n$.

\vskip 2mm
{\bf Subcase 2.3.1.3}. $z'\in \bigcup_{i=1}^{n-3}V(P_i)$.

W.l.o.g., assume that $z'\in V(P_{n-3})$. Recall that $z_n=\widetilde{z}$.

Let $F_1=EA_n^{n-1}\backslash \{z_1, \cdots, z_{n-4}, z_{n-2}, z_n, v\}$. According to Lemma \ref{lemEAn31}, $F_1$ is still connected. There is a $\{x_{n-1}', z_{n-3}, z\}$-tree $T_{n-1}'$ in $F_1$ and a $\{x_{n-2}', v', w\}$-tree $T_{n-2}'$ in $EA_n^{n-2}\backslash W$. Moreover, there is a $(z_{n-3}', w_{n-3}')$-path $\widetilde{P}$ in $EA_{n}^{n-3}$.

\begin{figure}[htbp]
\begin{minipage}[t]{0.5\linewidth}
\centering
\resizebox{0.95\textwidth}{!} {\includegraphics{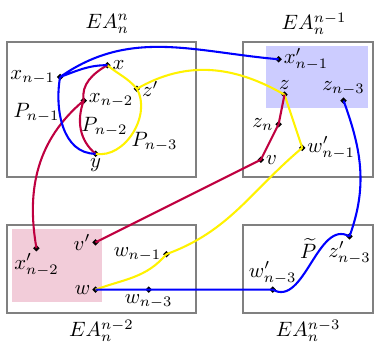}}
\caption{\small Illustration for Subcase 2.3.1.3} \label{figS211C2313}
\end{minipage}
\begin{minipage}[t]{0.5\linewidth}
\centering
\resizebox{0.95\textwidth}{!} {\includegraphics{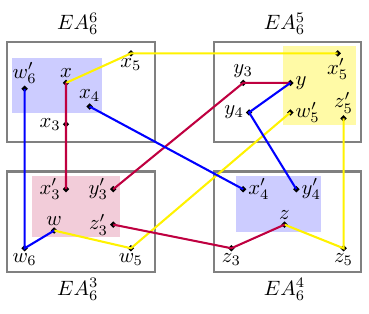}}
\caption{\small Illustration for Subcase 1.1.1} \label{figS1C111}
\end{minipage}
\end{figure}

Let
$$\widehat{T}_{n-3}= P_{n-3}\cup P[w,w_{n-1}']\cup \{zz', zw_{n-1}'\},$$
$$T_{n-2}=T_{n-2}'\cup P_{n-2}\cup \{x_{n-2}x_{n-2}', zz_n, z_nv, vv'\}$$
and
$$T_{n-1}=T_{n-1}' \cup P_{n-1}\cup \widetilde{P} \cup P[w,w_{n-3}']\cup \{x_{n-1}x_{n-1}', z_{n-3}z_{n-3}'\}.$$
See Figure \ref{figS211C2313}. Then $T_1, \cdots, T_{n-4},$ $\widehat{T}_{n-3}, T_{n-2}, T_{n-1}$ are ($n-1$) IDSTs in $EA_n$.

\vskip 2mm
{\bf Subcase 2.3.2}. $z'\notin \bigcup_{i=1}^{n-1}V(P_i)$.

Let $u_1, \cdots, u_{n-1}$ be $(n-1)$ neighbours of $z'$ in $EA_n^n$.

Firstly, we consider the case that $u_1\in \bigcup_{i=1}^{n-1}V(P_i)$. By similar methods in Subcase 2.3.1, we may construct $(n-1)$ IDSTs in $EA_n$ just by replacing the edge $zz'$ by the path $\{zz', z'u_1\}$ in $T_i$ when $u_1\in V(P_i)$, $n-3\le i\le n-1$.

Then we move to the case that $\{u_1, \cdots, u_{n-1}\}\cap \big(\bigcup_{i=1}^{n-1}V(P_i)\big)=\emptyset$. By Lemma \ref{lemEAn3}(2), let $u_{n-2}$ be a neighbour of $z'$ in $EA_n^n$ that $u_{n-2}'\in V(EA_n^{n-2})$. There is a $\{x_{n-1}', z, w_{n-1}'\}$-tree $T_{n-1}'$ in $EA_{n}^{n-1}\backslash \{z_1, \cdots, z_{n-3}\}$ and a $\{x_{n-2}', u_{n-2}', w\}$-tree $T_{n-2}'$ in $EA_n^{n-2}\backslash W$, respectively.

Let $$T_{n-1}=T_{n-1}'\cup P_{n-1}\cup P[w,w_{n-1}']\cup \{x_{n-1}x_{n-1}'\}$$
and
$$T_{n-2}=T_{n-2}'\cup P_{n-2}\cup \{x_{n-2}x_{n-2}', zz', z'u_{n-2}, u_{n-2}u_{n-2}'\}.$$
Then $T_1, \cdots,  T_{n-1}$ are ($n-1$) IDSTs in $EA_n$. \hfill$\Box$

\subsection{$\max\limits_{i\in [n]}\{|S\cap V(EA_n^{m:i})|\}=1$, $4\le m\le n$.}

\begin{lemma}\label{lemS1111}
For $n\ge 4$ and $4\le m\le n$, let $EA_n=EA_n^{m:1}\oplus EA_n^{m:2}\oplus \cdots \oplus EA_n^{m:n}$ and $S=\{x,y,z,w\}$ be any 4-subset of $V(EA_n)$. If  $|S\cap V(EA_n^{m:i})|\le 1$ for each $i\in[n]$. Then there are ($n-1$) IDSTs in $EA_n$.
\end{lemma}

\noindent{\bf Proof}\; For simplicity, we assume that $x=12\cdots n\in V(EA_n^n)$. We need to consider whether $\{y_n', z_n', w_n'\}\cap N_{EA_n^n}(x)=\emptyset$ or not.

\vskip 2mm

{\bf Case 1}. $y_n'\in N_{EA_n^n}(x)$.

\vskip 2mm
Since $|S\cap V(EA_n^{m:i})|\le 1$ for $4\le m\le n$ and each $i\in [n]$, Claim 1-Claim 3 may be inferred from Definition \ref{defEAn}.

\vskip 2mm

{\bf Claim 1}. $n\le 6$.

{\bf Proof of Claim 1}. Suppose to contrary that $n\ge 7$. We will conclude that there exists an integer $j$ ($4\le j\le n-1$) such that the $j$th symbol of $y$ is $j$.

For $y_n'\in \{x\circ(12), x\circ(123), x\circ(132)\}$, the $j$th symbol of $y_n$ is $j$ for each $4\le j\le n-1$ due to Definition \ref{defEAn}. Then there exists an integer $j$ ($4\le j\le n-1$) such that the $j$th symbol of $y$ is $j$ no matter $y_n\in \{y\circ(12), y\circ(123), y\circ(132)\}$ or $y_n\in \{y\circ(12)(3i)|4\le i\le n-1\}$.

For $y_n'\in \{x\circ(12)(3i)|4\le i\le n-1\}$, the $j$th symbol of $y_n$ is $j$ for each $j\in [n-1]\backslash \{1,2,3,i\}$. Then there always exists an integer $j$ ($4\le j\le n-1$) such that the $j$th symbol of $y$ is $j$ when $n\ge 7$.

The above analysis means that $\{x,y\}\subseteq V(EA_n^{j:j})$, which is a contradiction to the condition of this Lemma \ref{lemS1111}.  \hfill$\Box$

\vskip 2mm

{\bf Claim 2}. If $n=6$. Then $y\in\{215364, 214635\}$.

{\bf Proof of Claim 2}. By similar analysis in the proof of Claim 1, it has $y_n'\notin \{x\circ(12), x\circ(123), x\circ(132)\}$. Moreover, $y_n\notin \{y\circ(12), y\circ(123), y\circ(132)\}$. That means either $y_n'=214356$ or $y_n'=215436$. Then $y_n\in \{126354, 126435\}$ and $y\in\{215364, 214635\}$.  \hfill$\Box$

\vskip 2mm

{\bf Claim 3}. If $n= 5$. Then $y\in\{21453, 31452, 23451, 21534, 51234, 25134, 21354\}$.

{\bf Proof of Claim 3}. By similar above analysis, $y_n\notin \{y\circ(12), y\circ(123), y\circ(132)\}$ when $y_n'\in \{x\circ(12), x\circ(123), x\circ(132)\}$. Otherwise, $\{x,y\}\subseteq V(EA_5^{4:4})$, a contradiction to the assumption of this lemma. Therefore, we have $y=21453$ if $y_n'=x\circ(12)$, $y=31452$ if $y_n'=x\circ(132)$ and $y=23451$ if $y_n'=x\circ(123)$.

For $y_n'=x\circ(12)(34)=21435$. Then $y_n=12534$. Furthermore, we have $y\in \{21534, 51234, 25134, 21354\}$. \hfill$\Box$

\vskip 2mm
Based on Claim 1, we only need to focus on the case that $4\le n\le 6$.

\vskip 2mm

{\bf Subcase 1.1}. $n=6$.

According to Claim 2, $|\{z_n', w_n'\}\cap N_{EA_n^n}(x)|\le 1$.

\vskip 2mm

{\bf Subcase 1.1.1}. Either $z_n'\in N_{EA_n^n}(x)$ or $w_n'\in N_{EA_n^n}(x)$.

By symmetry, assume that $z_n'\in N_{EA_n^n}(x)$. Then $y=214635$ and $z=215364$ by Claim 2.

Note that $x_5=y_6'=215436$, $y_4'=215634\in V(EA_6^4)$, $x_4=z_6'$ and $z_5=125364$. Clearly, $y_4'\ne z_5$. Moreover, $w_6'\notin N_{EA_6^6}(x)$. By Table \ref{tablee1}, it is not difficult to check that, for any $w\in V(EA_n^i)$ ($1\le i\le 3$), $w$ cannot be adjacent to $y_i'$ and $z_i'$ at the same time.

\begin{table}[h] \footnotesize 
\centering
 \caption{\small }\label{tablee1}
\begin{tabular}{ccccccc}
\specialrule{0.05em}{3pt}{3pt}
 vertex $a$ & $a_1$  & $a_1'$ & $a_2$  & $a_2'$ & $a_3$  & $a_3'$ \\
\specialrule{0.05em}{3pt}{3pt}
$y=214635$ &  421635 & 245631   & 142635   & 415632  & 123645 & 215643 \\
\specialrule{0.05em}{3pt}{3pt}
$z=215364$ & 521364 &254361 & 152364 & 514362 & 123564& 214563\\
\specialrule{0.05em}{3pt}{3pt}
\end{tabular}
\end{table}

W.l.o.g., assume that $w\in V(EA_6^3)$ and $wy_3'\notin E(EA_6^3)$. That means $w_5'\ne y_3$.

Let $F_6=EA_6^6\backslash\{x_1, x_2, x_3, x_5\}$, $F_5=EA_6^5\backslash\{y_1, y_2, y_3, y_4\}$, $F_4=EA_6^4\backslash\{z_1, z_2, z_3,$ $z_5\}$ and $F_3=EA_6^3\backslash\{w_1, w_2, w_5, w_6\}$.

For $3\le i\le 6$, $F_i$ is connected. Hence, there is a $\{x_3',y_3',z_3',w\}$-tree $T_3'$ in $F_3$, a $\{x_4',y_4',z\}$-tree $T_4'$ in $F_4$, a $\{x_5',y, z_5',w_5'\}$-tree $T_5'$ in $F_5$ and a $\{x,x_4,w_6'\}$-tree $T_6'$ in $F_6$, respectively.

Let $T_3=T_3'\cup P[x,x_3']\cup P[y,y_3']\cup P[z,z_3']$, $T_4=T_4'\cup T_6'\cup  P[y,y_4']\cup P[w,w_6']\cup \{x_4x_4'\}$ and $T_5=T_5'\cup P[x,x_5']\cup P[z,z_5']\cup P[w,w_5']$. See Figure \ref{figS1C111}.

For $1\le i\le 2$, there is a $\{x_i',y_i',z_i',w_i'\}$-tree $T_i'$ in $EA_6^i$, hereafter, let $T_{i}=T_{i}'\cup P[x,x_{i}']\cup P[y,y_{i}']\cup P[z,z_{i}']\cup P[w,w_i'].$ Then $T_1, \cdots, T_5$ are five IDSTs in $EA_6$.

\vskip 2mm

{\bf Subcase 1.1.2}. $z_n'\notin N_{EA_n^n}(x)$ and $w_n'\notin N_{EA_n^n}(x)$.

W.l.o.g., assume that $y\in V(EA_6^5)$, $z\in V(EA_6^4)$ and $w\in V(EA_6^3)$.

Then $z_6'\ne x_4$. We only need to consider the case that $z_5'\ne y_4$ and $w_5'\ne y_3$. Otherwise, the proof is similar to that of Subcase 1.1.1.

Under this circumstance, we may obtain five IDSTs $T_1, \cdots, T_5$ again by the same method in Subcase 1.1.1.

\vskip 2mm

{\bf Subcase 1.2}. $n=5$.

By Claim 3, it has $|\{y_5', z_5', w_5'\}\cap N_{EA_5^5}(x)|\le 2$. Otherwise, either $|\{y,z,w\}\cap V(EA_5^{4:3})|\ge 2$ or $|\{y,z,w\}\cap V(EA_5^{4:5})|\ge 2$. Hence, assume that $w_5'\notin N_{EA_5^5}(x)$.

\vskip 2mm

{\bf Subcase 1.2.1}. $z_5'\in N_{EA_5^5}(x)$.

By Claim 3 and the fact that $|S\cap V(EA_5^{m:i})|\le 1$ for $4\le m\le 5$ and $i\in [5]$, it has $y\in\{{21453, 31452, 23451}\}$ and $z\in\{{21534, 51234, 25134}\}$. Note that $\{y',z\}\subseteq V(EA_5^4)$. Moreover, it has $z_5=x_4'=12534$.

\vskip 2mm

{\bf Subcase 1.2.1.1}. $y=21453$.

W.l.o.g., assume that $w\in V(EA_5^2)$. Note that $x'\in V(EA_5^3)$ and $y_5=x_3'=12543$. Furthermore, $z_3\ne y_4'$ when $z\in\{21534,51234,25134\}$.

Let $F_5=EA_5^5\backslash \{x_1,x_2,x_3\}$, $F_4=EA_5^4\backslash \{z_1,z_2,z_3\}$, $F_3=EA_5^3\backslash \{y_1,y_4,y_5\}$ and $F_2=EA_5^2\backslash \{w_1,w_3,w_5\}$.

Clearly, $F_i$ is connected for $2\le i\le 5$. Thus, there is a $\{x,x_4,w_5'\}$-tree $T_5'$ in $F_5$, a $\{x_4', y_4', z\}$-tree $T_4'$ in $F_4$, a $\{x', y, z_3', w_3'\}$-tree $T_3'$ in $F_3$ and a $\{x_2', z_2', w\}$-tree $T_2'$ in $F_2$, respectively.

\begin{figure}[htbp]
\begin{minipage}[t]{0.5\linewidth}
\centering
\resizebox{0.95\textwidth}{!} {\includegraphics{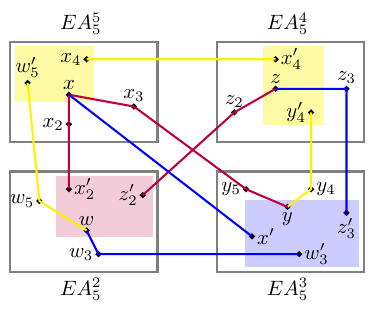}}
\caption{\small Illustration for Subcase 1.2.1.1} \label{figS1C1211}
\end{minipage}
\begin{minipage}[t]{0.5\linewidth}
\centering
\resizebox{0.95\textwidth}{!} {\includegraphics{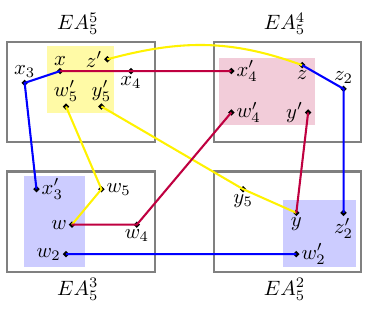}}
\caption{\small Illustration for Subcase 1.2.1.2.1} \label{figS1C12121}
\end{minipage}
\end{figure}

Let $T_2=T_2'\cup P[x,x_2']\cup P[x,x_3']\cup P[z,z_2']\cup \{x_3'y\}$, $T_3=T_3'\cup P[z,z_3']\cup P[w,w_3']\cup \{xx'\}$ and $T_4=T_4'\cup T_5'\cup P[y,y_4'] \cup P[w,w_5']\cup \{x_4x_4'\}$. See Figure \ref{figS1C1211}. Let $T_1'$ and $T_1$ be the same as in Subcase 1.1.1. Then $T_1, \cdots, T_4$ are four desired IDSTs in $EA_5$.

{\bf Subcase 1.2.1.2}. $y=31452$.

\begin{table}[htbp] \footnotesize
\begin{minipage}[t]{0.5\linewidth}
\centering
 \caption{\small }\label{table1}
\begin{tabular}{ccccc}
\specialrule{0.05em}{3pt}{3pt}
 vertex $a$ & $a_1$  & $a_1'$  & $a_3$  & $a_3'$ \\
\specialrule{0.05em}{3pt}{3pt}
$y=31452$ & 43152  &  34251 &  14352 & 41253 \\
\specialrule{0.05em}{3pt}{3pt}
$z=25134$ &52134 & 25431 & 52314  &25413 \\
\specialrule{0.05em}{3pt}{3pt}
$z=51234$ &25134 & 52431 & 15324  & 51423\\
\specialrule{0.05em}{3pt}{3pt}
$z=21534$&52134 & 25431 & 12354  & 21453 \\
\specialrule{0.05em}{3pt}{3pt}
\end{tabular}

\end{minipage}
\quad
\begin{minipage}[t]{0.5\linewidth}
\centering
 \caption{\small }\label{table2}
\begin{tabular}{ccccc}
\specialrule{0.05em}{3pt}{3pt}
 vertex $a$ & $a_2$  & $a_2'$  &  $a_3$  & $a_3'$ \\
\specialrule{0.05em}{3pt}{3pt}
$y=23451$ & 34251  &  43152 &  42351 & 24153 \\
\specialrule{0.05em}{3pt}{3pt}
$z=25134$ & 51234  & 15432  &  52314  & 25413 \\
\specialrule{0.05em}{3pt}{3pt}
$z=51234$ & 15234  &  51432 & 15324   & 51423 \\
\specialrule{0.05em}{3pt}{3pt}
$z=21534$& 15234 & 51432  &  12354  & 21453  \\
\specialrule{0.05em}{3pt}{3pt}
\end{tabular}

\end{minipage}
\end{table}

It is seen from Table \ref{table1} that the 3rd and 4th symbols of $y_1'$ are different from those of $z_1'$ no matter $z=25134$ or $z=51234$ or  $z=21534$. By Definition \ref{defEAn}, for any vertex $u\in V(EA_5^1)$, it cannot be adjacent to $y_1'$ and $z_1'$ at the same time.

When $z=25134$, note that $z_3'=25413$ and $y_3'=41253$. For any vertex $u\in V(EA_5^3)$, it is has $|\{y_3', z_3'\}\cap N_{EA_5^3}(u)|\le 1$. When $z=51234$, a vertex $u\in V(EA_5^3)$ may be adjacent to $y_3'$ and $z_3'$ at the same time if and only if $u=14523$.

\vskip 2mm

{\bf Subcase 1.2.1.2.1}. $w\in V(EA_5^i)$ is a common neighbour of $y_i'$ and $z_i'$, $i\in \{1,3\}$.

By above analysis, there are only two possibilities: (1) $z=51234$ and $w=14523$; (2) $z=21534$ and $w\in V(EA_5^3)$.

Note that $\{x,y,z,w\}\in V(AN_5^1)$ when $z=51234$ and $w=14523$. Then there are four IDSTs in $EA_5$ by Lemma \ref{lemS4}.

Now we consider the possibility that $z=21534$  and $w\in V(EA_5^3)$. Remind that $w_5'\ne x_3$. Note that $z'\in V(EA_5^5)$, $x_2'=y_5$, $x_4'=z_5$ and $z_3'=w_4$.

Let $F_5=EA_5^5\backslash \{x_1, x_3, x_4\}$, $F_4=EA_5^4\backslash \{z_1, z_2\}$, $F_3=EA_5^3\backslash \{w_1, w_4, w_5\}$ and $F_2=EA_5^2\backslash \{y_1, y_5\}$. Note that $F_i$ is connected for $2\le i\le 5$. Thus, there is a $\{x, y_5', z', w_5'\}$-tree $T_5'$ in $F_5$, a $\{x_4', y', z, w_4'\}$-tree $T_4'$ in $F_4$, a $\{x_3', w, w_2\}$-tree $T_3'$ in $F_3$ and a $\{y, z_2', w_2'\}$-tree $T_2'$ in $F_2$.

Let $T_2=T_2'\cup T_3'\cup P[x,x_3']\cup P[z,z_2']\cup \{w_2w_2'\}$, $T_4=T_4'\cup P[x,x_4']\cup P[w, w_4']\ \{yy'\}$ and $T_5=T_5'\cup P[y, y_5']\cup P[w, w_5']\cup \{zz'\}$. See Figure \ref{figS1C12121}. Let $T_1$ be the same as in Subcase 1.1.1. Then $T_1, T_2, T_4$ and $T_5$ are four IDSTs in $EA_5$.

\vskip 2mm

{\bf Subcase 1.2.1.2.2}. Either $wy_i'\notin E(EA_5^i)$ or $wz_i'\notin E(EA_5^i)$, where $w\in V(EA_5^i)$ and $i\in \{1,3\}$.

W.l.o.g., assume that $w\in V(EA_5^3)$ and $wy_3'\notin E(EA_5^3)$. That means $y_3'\ne w_2$.

Remind that $x_3'\ne w_5$, $z\in V(EA_5^4)$ and $x_4=z_5'$. Moreover, it has $y'\ne z_2$ when $z\in\{51234,25134,21534\}$.

There is a $\{x, z_5', w_5'\}$-tree $T_5'$ in $EA_5^5\backslash \{x_1, x_2, x_3\}$, a $\{y', z, z_5\}$-tree $T_4'$ in $EA_5^4\backslash \{z_1,$ $z_2, z_3\}$, a $\{x_3', y_3', z_3', w\}$-tree $T_3'$ in $EA_5^3\backslash \{w_1, w_2, w_5\}$ and a $\{x_2', y, z_2', w_2'\}$-tree $T_2'$ in $EA_5^2\backslash \{y_1, y_3\}$, respectively.

\begin{figure}[htbp]
\begin{minipage}[t]{0.5\linewidth}
\centering
\resizebox{0.92\textwidth}{!} {\includegraphics{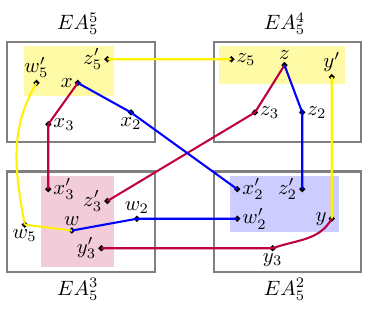}}
\caption{\small Illustration for Subcase 1.2.1.2.2} \label{figS1C12122}
\end{minipage}
\begin{minipage}[t]{0.5\linewidth}
\centering
\resizebox{0.99\textwidth}{!} {\includegraphics{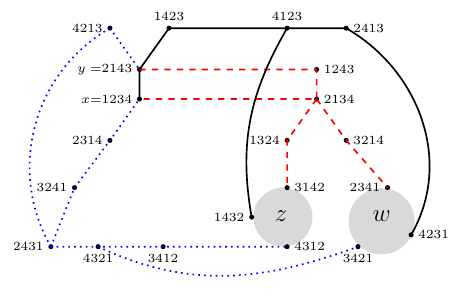}}
\caption{\small Three IDSTs in $EA_4$} \label{figEA4zR}
\end{minipage}
\end{figure}

Let $T_4=T_4'\cup T_5'\cup P[w, w_5']\cup \{yy', z_5z_5'\}$, $T_3=T_3'\cup P[x, x_3']\cup P[y, y_3']\cup P[z, z_3']$ and $T_2=T_2'\cup P[x, x_2']\cup P[z,z_2']\cup P[w, w_2']$.
See Figure \ref{figS1C12122}. Let $T_1$ be the same as in Subcase 1.1.1. Then $T_1, T_2, T_3$ and $T_4$ are four IDSTs in $EA_5$.

\vskip 2mm

{\bf Subcase 1.2.1.3}. $y=23451$.

Based on Table \ref{table2} and similar arguments in Subcase 1.2.1.2, for $2\le i\le 3$ and any vertex $w\in V(EA_5^i)$, $w$ is adjacent to both $y_i'$ and $z_i'$ if and only if one of the two possibilities: (1) $z=25134$ and $w=42513$; (2) $z=21534$ and $w\in V(EA_5^3)$.

Note that $\{x, y, z, w\}\subseteq V(AN_5^1)$ when $z=25134$ and $w=42513$, then there are four IDSTs in $EA_5$ according to Lemma \ref{lemS4}. Hence, we only need to consider the possibility that $z=21534$ and $w\in V(EA_5^3)$ if $w\in V(EA_5^i)$ is a common neighbour of $y_i'$ and $z_i'$, $2\le i\le 3$. Note that $y'\in V(EA_5^4)$ and $y'\ne z_1$ when $z=21534$, moreover, $z'\in V(EA_5^5)$. By similar analysis in Subcase 1.2.1.2.1, four IDSTs can be obtained in $EA_5$ and the proof is omitted.

If $wy_i'\notin E(EA_5^i)$ or $wz_i'\notin E(EA_5^i)$, where $w\in V(EA_5^i)$ and $i\in \{2,3\}$. Then we may obtain four IDSTs by similar discussions in Subcase 1.2.1.2.2.

\vskip 2mm

{\bf Subcase 1.2.2}. $z_5'\notin N_{EA_5^5}(x)$.

W.l.o.g., assume that $z\in V(EA_5^4)$, $w\in V(EA_5^3)$ and $y\in V(EA_5^2)$. Moreover, we may assume that $\{z_2', w_2'\}\cap N_{EA_5^2}(y)=\emptyset$. Otherwise, four IDSTs can be obtained by similar method in Subcase 1.2.1.

Note that there is a $\{x_2', y, z_2', w_2'\}$-tree $T_2'$ in $EA_5^2\backslash \{y_1, y_3, y_4\}$ and a $\{y_4', z, z_5\}$-tree $T_4'$ in $EA_5^4\backslash \{z_1, z_2, z_3\}$. Let $T_i'$ ($i\in\{1,3,5\}$) and $T_i$ ($i\in [3]$) be the same as in Subcase 1.2.1.2.2. Set $T_4=T_4'\cup T_5'\cup P[w, w_5']\cup P[y, y_4']\cup \{z_5z_5'\}$. Then $T_1, \cdots, T_4$ are four IDSTs in $EA_5$.

\vskip 2mm

{\bf Subcase 1.3}. $n=4$.

W.l.o.g., assume that $y\in V(EA_4^3)$, $z\in V(EA_4^2)$ and $w\in V(EA_4^1)$. Then $y\in \{2143, 4123, 2413\}$ by Figure \ref{figEA4}. By symmetry, we only need to consider the case that $y=2143$ or $y=2413$.

\vskip 2mm

{\bf Subcase 1.3.1}. $y=2143$.

Recall that $EA_4=AN_4^1\otimes AN_4^2$ and $\{x,y\}\subseteq V(AN_4^1)$. We conclude from Lemma \ref{lemS4} that three IDSTs may be contructed if $\{z,w\}\subseteq V(AN_4^1)$. Moreover, if $\{z, w\}\cap V(AN_4^1)=\{z\}$, then three IDSTs may be obtained by Lemma \ref{lemANS3} since $\widetilde{w}\in V(EA_4^1)$ and $|E(\widetilde{w}, \{x,y,z\})|\le 1$. Therefore, we only need to consider the case that $\{z,w\}\subseteq V(AN_4^2)$.

Let $V_1=\{1432,4312,3142\}$ and $V_2=\{2341, 3421, 4231\}$. Then $z\in V_1$ and $w\in V_2$.

Note that the subgraph induced on $V_i$ is a clique, $1\le i\le 2$. By viewing $V_i$ ($1\le i\le 2$) as a vertex, Figure \ref{figEA4zR} shows that there are three internally edge disjoint $\{x,y,V_1,V_2\}$-trees in $EA_4$, where three trees are represented by black solid lines, blue dotted lines and red dashed lines, respectively.

From the trees shown in Figure \ref{figEA4zR}, we can obtain three desired IDSTs in $EA_4$ by adding two edges $z(z\circ(123))$ and $z(z\circ(132))$ in $V_1$ and by adding two edges $w(w\circ(123))$ and $w(w\circ(132))$ in $V_2$.

\begin{figure}[htbp]
\begin{minipage}[t]{0.49\linewidth}
\centering
\resizebox{0.98\textwidth}{!} {\includegraphics{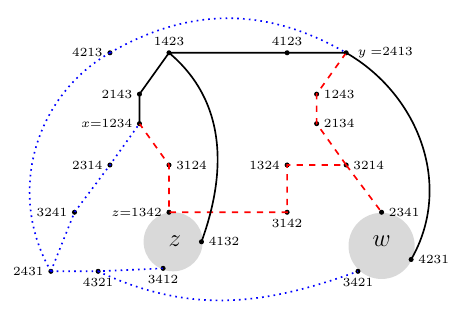}}
\caption{\small Three IDSTs in $EA_4$} \label{figEA4zL}
\end{minipage}
\begin{minipage}[t]{0.49\linewidth}
\centering
\resizebox{0.95\textwidth}{!} {\includegraphics{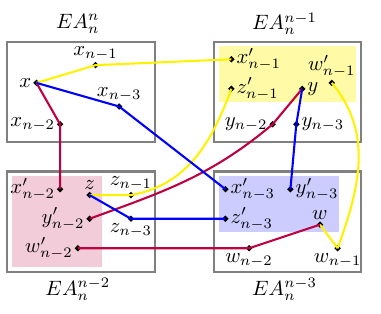}}
\caption{\small Illustration for Case 2} \label{figS1C2}
\end{minipage}
\end{figure}

{\bf Subcase 1.3.2}. $y=2413$.

By similar analysis in Subcase 1.3.1, we only need to consider the possibility that $z\in V(AN_4^1)$ and $w\in V(AN_4^2)$. Figure \ref{figEA4zL} shows that we can find three IDSTs in $EA_4$ by adding four corresponding edges when $z\in\{1342,4132,3412\}$ and $w\in\{2341, 3421, 4231\}$.

\vskip 2mm

{\bf Case 2}. $\{y_n', z_n', w_n'\}\cap N_{EA_n^n}(x)=\emptyset$.

W.l.o.g., assume that $y\in V(EA_n^{n-1})$, $z\in V(EA_n^{n-2})$ and $w\in V(EA_n^{n-3})$. Furthermore, we may assume that $\{z_{n-1}', w_{n-1}'\}\cap N_{EA_n^{n-1}}(y)=\emptyset$ and $z_{n-3}\ne w_{n-2}'$. Otherwise, similar arguments in Case 1 may be applied to obtain ($n-1$) IDSTs in $EA_n$.

For $i\in[n-4]$, there is a $\{x_i', y_i', z_i', w_i'\}$-tree $T_{i}'$ in $EA_n^i$ since $\{x_i', y_i', z_i', w_i'\}\subseteq V(EA_n^i)$. Let
$$T_{i}=T_{i}'\cup P[x,x_{i}']\cup P[y,y_{i}']\cup P[z,z_{i}']\cup P[w,w_i'], ~~~~ i\in[n-4].$$

Denote
$$F_{n-1}=EA_n^{n-1}\backslash \{y_1,\cdots, y_{n-2}\}, \quad F_{n-2}=EA_n^{n-2}\backslash\{z_1,\cdots, z_{n-3}, z_{n-1}\}$$
and $F_{n-3}=EA_n^{n-3}\backslash\{w_1,\cdots, w_{n-4}, w_{n-2}, w_{n-1}\}.$

Clearly, the subgraph $F_i$ is connected for $n-3\le i\le n-1$. Hereafter, there is a $\{x_{n-1}', y, z_{n-1}', w_{n-1}'\}$-tree $T_{n-1}'$ in $F_{n-1}$, a $\{x_{n-2}', y_{n-2}', z, w_{n-2}'\}$-tree $T_{n-2}'$ in $F_{n-2}$ and a $\{x_{n-3}', y_{n-3}', z_{n-3}', w\}$-tree $T_{n-3}'$ in $F_{n-3}$, respectively.

Let
$$T_{n-1}=T_{n-1}'\cup P[x,x_{n-1}']\cup P[z,z_{n-1}']\cup P[w,w_{n-1}'],$$
$$T_{n-2}=T_{n-2}'\cup P[x,x_{n-2}']\cup P[y,y_{n-2}']\cup P[w,w_{n-2}']$$
and
$$T_{n-3}=T_{n-3}'\cup P[x,x_{n-3}']\cup P[y,y_{n-3}']\cup P[z,z_{n-3}'].$$
See Figure \ref{figS1C2}. Then $T_1,\cdots, T_{n-1}$ are ($n-1$) desired IDSTs in $EA_n$. \hfill$\Box$

\section{Proof of Theorem \ref{thk4EAn}}\label{secK4EAn}

In this section, we are going to prove Theorem \ref{thk4EAn} by induction on $n$. Firstly, we focus on the induction basis for $n=3$.

\begin{lemma}\label{lemEA3}
 $\kappa_4(EA_3)\ge 2$.
\end{lemma}

\noindent{\bf Proof}\; It suffices to prove that there are two IDSTs in $EA_3$ for any 4-subset $S=\{x,y,z,w\}\subseteq V(EA_3)$. Note that $EA_3=AN_3^1\otimes AN_3^2$ and $|V(AN_3^i)|=3$ for $1\le i\le 2$. The following two cases are discussed.

{\bf Case 1.} $|S\cap V(AN_3^i)|=3$ for $1\le i\le 2$.

\begin{figure}[htbp]
\begin{minipage}[b]{0.31\linewidth}
\centerline {\includegraphics[width=\textwidth]{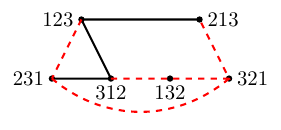}}
\caption{\small Two IDSTs\\ in $EA_3$} \label{figEA31}
\end{minipage}
\begin{minipage}[b]{0.31\linewidth}
\centering
\centerline {\includegraphics[width=\textwidth]{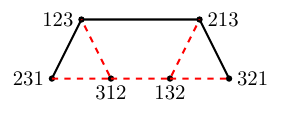}}
\caption{\small Two IDSTs\\ in $EA_3$} \label{figEA321}
\end{minipage}
\begin{minipage}[b]{0.31\linewidth}
\centering
\centerline {\includegraphics[width=\textwidth]{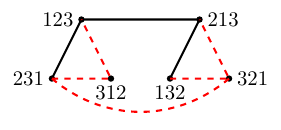}}
\caption{\small Two IDSTs\\ in $EA_3$} \label{figEA322}
\end{minipage}
\end{figure}

By symmetry, we may assume that $x=123$, $y=231$, $z=312$ and $w=213$. Then two IDSTs are shown in Figure \ref{figEA31}, where two trees are represented by black solid lines and red dashed lines, respectively.

{\bf Case 2.} $|S\cap V(AN_3^1)|=|S\cap V(AN_3^2)|=2$.

W.l.o.g., we may assume that $x=123$ and $y=231$. By investigating the structure of $EA_3$, two IDSTs are shown in Figure \ref{figEA321} if $\{z,w\}=\{213,321\}$ and in Figure \ref{figEA322} if $\{z,w\}=\{213,132\}$, respectively.   \hfill$\Box$

\vskip 2mm

Now, we are in a position to prove our main result.
\vskip 2mm

\noindent{\bf Proof of Theorem \ref{thk4EAn}}. Combined with Lemma \ref{lemupperKk} and Lemma \ref{lemEAn1}(1), $\kappa_4(EA_n)\le \delta(EA_n)-1=n-1$ for $n\ge 3$. We shall prove the reverse inequality by induction on $n$. First of all, $\kappa_4(EA_3)\ge 2$ by Lemma \ref{lemEA3}. Now suppose that $n\ge 4$ and the result is true for any integer $l<n$, i.e., $\kappa_4(EA_l)\ge l-1$. Let $S=\{x,y,z,w\}$ be any 4-subset of $V(EA_n)$.

Remind that $EA_n=EA_n^{m:1}\oplus EA_n^{m:2}\oplus \cdots \oplus EA_n^{m:n}$, $4\le m\le n$. The following cases are distinguished.

\vskip 2mm
{\bf Case 1.} There exists an integer $i\in [n]$ such that $S\subseteq V(EA_n^{m:i})$.

W.l.o.g., we may assume that $S\subseteq V(EA_n^n)$. According to the fact that $EA_n^n$ is isomorphic to $EA_{n-1}$ and our induction hypothesis, there are ($n-2$) IDSTs $T_1, \cdots, T_{n-2}$ in $EA_n^n$.

By Lemma \ref{lemEAn4}, there exists a $\{x', y', z', w'\}$-tree $T_{n-1}'$ in $\bigcup_{i=1}^{n-1}EA_n^i$. Let $T_{n-1}=T_{n-1}'\cup\{xx', yy', zz', ww'\}$. Clearly, $T_1, \cdots, T_{n-2}, T_{n-1}$ are ($n-1$) IDSTs in $EA_n$.

\vskip 2mm
{\bf Case 2.} There exists an integer $i\in [n]$ such that $|S\cap V(EA_n^{m:i})|=3$.

With the help of Lemma \ref{lemS3}, we may get ($n-1$) IDSTs in $EA_n$.

\vskip 2mm
{\bf Case 3.} There exist different integers $i$ and $j$ in $[n]$ such that $|S\cap V(EA_n^{m:i})|=|S\cap V(EA_n^{m:j})|=2$.

Then ($n-1$) IDSTs in $EA_n$ can be constructed due to Lemma \ref{lemS22}.

\vskip 2mm
{\bf Case 4.} There exist different integers $i$, $j$ and $k$ in $[n]$ such that $|S\cap V(EA_n^{m:i})|=2$ and $|S\cap V(EA_n^{m:j})|=|S\cap V(EA_n^{m:k})|=1$.

Based on Lemma \ref{lemS211}, there are ($n-1$) IDSTs in $EA_n$.

\vskip 2mm
{\bf Case 5.} $|S\cap V(EA_n^{m:i})|\le 1$ for each $i\in[n]$.

By Lemma \ref{lemS1111}, ($n-1$) IDSTs in $EA_n$ can be obtained. The proof is done. \hfill$\Box$

\section{Conclusion}\label{seccon}

The generalized $k$-connectivity is a natural generalization of the classical connectivity and can serve for measuring the capability of a network $G$ to connect any $k$ vertices in $G$. In this paper, we focus on the generalized 4-connectivity of the $n$-dimensional godan graph $EA_n$. Together with Theorem \ref{thk4EAn} and Lemma \ref{lemKkk-1}, the generalized 3-connectivity of $EA_n$ can be derived immediately. This result was also investigated in \cite{BPN3}.

%
%

%



\end{document}